\documentclass{jams-l}

\usepackage{amssymb}
\usepackage{amsmath}

\usepackage{graphicx}

\usepackage[cmtip,all]{xy}

\usepackage[dvipsnames]{xcolor}
\usepackage{stmaryrd}
\usepackage{comment}
\usepackage{tikz}
\usepackage{caption}
\usepackage{subcaption}
\usepackage{float}
\usepackage{wrapfig}
\usepackage{bm}
\usepackage{dashbox}
\usepackage{multirow}
\usepackage{xcolor}
\usepackage{hyperref}
\usepackage{cancel}
\usepackage{mathrsfs}
\usepackage{fontawesome5}

\usepackage{todonotes}

\usetikzlibrary{quotes,angles}

\newtheorem{theorem}{Theorem}[section]
\newtheorem{lemma}[theorem]{Lemma}
\newtheorem{proposition}[theorem]{Proposition}
\newtheorem{proposition-definition}[theorem]{Proposition-Definition}
\newtheorem{corollary}[theorem]{Corollary}

\theoremstyle{definition}
\newtheorem{definition}[theorem]{Definition}

\newtheorem{assumption}[theorem]{Assumption}

\theoremstyle{remark}
\newtheorem{remark}[theorem]{Remark}

\numberwithin{equation}{section}

\DeclareMathOperator{\id}{id}

\DeclareMathOperator{\vol}{Vol}

\DeclareMathOperator{\Real}{Re}
\DeclareMathOperator{\Imag}{Im}

\DeclareMathOperator{\vect}{Vect}

\DeclareMathOperator{\covol}{Covol}

\DeclareMathOperator{\rank}{rank}

\newcommand{\transpose}[1]{\phantom{}^{t} #1}

\newcommand{\hooklongrightarrow}{\lhook\joinrel\longrightarrow}

\newcommand{\twoheadlongrightarrow}{\relbar\joinrel\twoheadrightarrow}

\begin{document}

\title[A Minkowski-type theorem on distances to cusps]{A Minkowski-type theorem on distances to cusps: the class number one case}


\author{Mathieu Dutour}
\address{Institut Camille Jordan}
\curraddr{}
\email{mathieu.dutour@univ-st-etienne.fr}
\thanks{}



\date{}

\dedicatory{}

\begin{abstract}
    In the study of Euclidean lattices, the product of the successive minima is bounded from above and below by explicit quantities. This result is known as Minkowski's second theorem, and can be refined to include Hermite's constant in the upper bound, which measures how short a non-zero vector can be in a given lattice. A version of this result exists in the context of number fields, where lattices are replaced with rigid adelic spaces, and successive minima with the Roy--Thunder minima. In this paper, drawing on the analogy between rank $2$ Euclidean lattices and points in $\mathbb{H}$, we will see an analogy between $2$-dimensional rigid adelic spaces and points in $\mathbb{H}^n$, and use that to translate the Minkowski-type theorem on Roy--Thunder minima into a theorem on the distances to cusps in $\mathbb{H}^n$.
\end{abstract}

\maketitle

\tableofcontents


\section{Introduction}

    \subsection{\texorpdfstring{Situation for rank $2$ lattices}{Situation for rank 2 lattices}}

        Consider a rank $2$ Euclidean lattice $E$, \textit{i.e.} a free $\mathbb{Z}$-module of rank $2$ included in $\mathbb{R}^2$, and denote by $\left \Vert \cdot \right \Vert$ the usual Euclidean norm. The first minimum of $E$ is defined as
        \begin{equation}
            \begin{array}{lll}
                \lambda_1 \left( E \right) & = & \min \left\{ \left \Vert v \right \Vert \; \middle \vert \; v \in E \setminus \left \{ 0 \right \} \right\},
            \end{array}
        \end{equation}
        \textit{i.e.} as the shortest length of a non-zero vector in $E$. We can use this first minimum of rank $2$ lattices to define Hermite's constant as
        \begin{equation}
            \begin{array}[t]{lll}
                \gamma_2 & = & \displaystyle \raisebox{4pt}{$\max\limits_{\substack{E \text{\, lattice}\phantom{{}^2} \\[0.2em] \rank E \, = \, 2}}$} \; \frac{\lambda_1 \left( E \right)^2}{\covol E},
            \end{array}
        \end{equation}
        and it can be shown that we have
        \begin{equation}
            \begin{array}{lll}
                \gamma_2 & = & \displaystyle \frac{2}{\sqrt{3}},
            \end{array}
        \end{equation}
        this maximum being reached by the hexagonal lattice. Instead of considering arbitrary (full-rank) lattices $E \subset \mathbb{R}^2$, one can also consider, without any loss of generality, the standard lattice $\mathbb{Z}^2$, endowed with a twist of the usual Euclidean inner-product $\left< \cdot, \cdot \right>$ by a symmetric positive-definite matrix $S \in S_2^{++} \left( \mathbb{R} \right)$. The twisted inner-product is then defined as
        \begin{equation}
        \label{eq:twistedInnerProduct}
            \begin{array}{lll}
                \displaystyle \left< \cdot, \cdot \right>_S & = & \displaystyle \left< S \cdot, \cdot \right>.
            \end{array}
        \end{equation}
        Such matrices $S$ can be identified with points in $\mathbb{H} \times \mathbb{R}_+^{\ast}$, using the classical transformations (see section \ref{subsec:posdeftoH} for details)
        \begin{equation}
        \label{eq:correspondanceS2H}
            \begin{array}{ccccccc}
                \multicolumn{3}{c}{S_2^{++} \left( \mathbb{R} \right)} & \longrightarrow & \multicolumn{3}{c}{\mathbb{H} \times \mathbb{R}_+^{\ast}} \\[1em]

                S & = & \left( \begin{array}{cc} u & v \\[0.4em] v & w \end{array} \right) & \longmapsto & \displaystyle \left( \tau_S, \det S \right) & = & \displaystyle \left( \frac{v + i \sqrt{\det S}}{w}, \, \det S \right) \\[2em]

                \multicolumn{3}{c}{\displaystyle \frac{\sqrt{\lambda}}{y} \left( \begin{array}{cc} x^2+y^2 & x \\[0.4em] x & 1 \end{array} \right)} & \longmapsfrom & \multicolumn{3}{c}{\displaystyle \left( \tau \, = \, x + i y, \, \lambda \right)}
            \end{array}
        \end{equation}
        and these transformations are compatible with the natural actions of $PSL_2 \left( \mathbb{Z} \right)$

        \begin{itemize}
            \item by congruence on $S_2^{++} \left( \mathbb{R} \right)$;
            \item by homographic transformations on $\mathbb{H}$, and the trivial action on $\mathbb{R}_+^{\ast}$.
        \end{itemize}
        For $S \in S_2^{++} \left( \mathbb{R} \right)$, denote by $E_S$ the lattice $\mathbb{Z}^2$ with inner-product as in \eqref{eq:twistedInnerProduct}. Since the first minimum $\lambda_1 \left( E_S \right)$ does not depend on the congruence class of $S$, one should be able to translate it into an invariant related to the class of the point $\tau_S$ in the modular curve $PSL_2 \left( \mathbb{Z} \right) \backslash \mathbb{H}$. A direct computation shows that we have
        \begin{equation}
        \label{eq:linkLambda1Mu}
            \begin{array}{lll}
                \lambda_1 \left( E_S \right)^2 & = & \left( \max\limits_{\gamma \in PSL_2 \left( \mathbb{Z} \right)} \Imag \left( \gamma \cdot \tau_S \right) \right)^{-1}.
            \end{array}
        \end{equation}
        Therefore, the first minimum of $E_S$ is closely related with the distance between $\tau_S$ and its closest cusp for the action of $PSL_2 \left( \mathbb{Z} \right)$. The relation \eqref{eq:linkLambda1Mu} has a geometric interpretation. On the usual fundamental domain for the action of $PSL_2 \left( \mathbb{Z} \right)$ on $\mathbb{H}$, which is the sphere of influence\footnote{See definition \ref{def:sphereInfluence}.} of the cusp $\infty$, represented in figure \ref{fig:fundamentalDomainPSL2Z} below,
        
        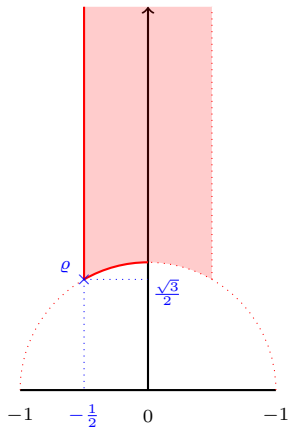
\begin{figure}[H]
    
            \centering
        
                \begin{tikzpicture}[scale=1.7]
        
                    \draw[thick] (-1,0) coordinate (1) node {} -- (1,0) coordinate (2) node {};
                    \draw[thick, ->] (0,0) -- (0,3);
                
                    \node (A) at (-0.5, 0.8660254) {};
                    \node (B) at (0.5,0.8660254) {};
                
                    \draw[thick,red] (0,1) arc (90:120:1cm);
                    \draw[dotted,red] (1,0) arc (0:90:1cm);
                
                    \draw[thick,red] (-0.5,0.8660254) -- (-0.5,3);
                    \draw[dotted,red] (0.5,0.8660254) -- (0.5,3);
                
                    \draw[dotted,red] (-0.5,0.8660254) arc (120:180:1cm);

                    \draw[dotted,blue] (-1/2,0.866025) -- (0,0.866025);
                    \draw[dotted,blue] (-1/2,0.866025) -- (-1/2,0);
                
                    \node (C) at (-1,-0.2) {$\scriptstyle -1$};
                    \node (D) at (1,-0.2) {$\scriptstyle -1$};
                    \node (E) at (0,-0.2) {$\scriptstyle 0$};
                    \node (E) at (-1/2,0.866025) {$\scriptstyle \color{blue} \times$};
                    \node (F) at (-1/2-0.15,0.866025+0.1) {$\scriptstyle \color{blue} \varrho$};
                    \node (G) at (0.15,0.866025-0.1) {$\scriptstyle \color{blue} \frac{\sqrt{3}}{2}$};
                    \node (H) at (-1/2,-0.2) {$\scriptstyle \color{blue} -\frac{1}{2}$};
                
                    \fill[red,opacity=0.2] (A) arc (120:60:1cm) -- (0.5,3) -- (-0.5,3) -- (A) -- cycle;

                \end{tikzpicture}
        
            \caption{Fundamental domain for the action of $PSL_2 \left( \mathbb{Z} \right)$ on $\mathbb{H}$}
            \label{fig:fundamentalDomainPSL2Z}
    
            \end{figure}
        
        \noindent the point which is the farthest away from $\infty$ is
        \begin{equation}
            \begin{array}{lll}
                \varrho & = & \left( - \frac{1}{2}, \; \frac{\sqrt{3}}{2} \right).
            \end{array}
        \end{equation}
        The inverse of its $y$-coordinate, which is the square of its distance to $\infty$, is the constant $\gamma_2$. We can also translate the second minimum of a rank $2$ lattice, defined~as
        \begin{equation}
            \begin{array}{lll}
                \lambda_2 \left( E \right) & = & \min \left \{ r>0 \; \middle \vert \; \dim \vect \left( \left \{ v \in E, \; \left \Vert v \right \Vert \, \leqslant \, r \right \} \right) \, = \, 2 \right \},
            \end{array}
        \end{equation}
        into an invariant related to the distance between~$\tau_S$ and its second closest cusp. Using these two interpretations of $\lambda_1 \left( E \right)$ and $\lambda_2 \left( E \right)$, it is possible to give a formulation of Minkowski's second theorem in terms of distances to cusps.

    \subsection{Statement of the results}

        The main purpose of this paper is to generalise the observations made in the previous section to the case of a totally real number field~$K$ of degree $n$, whose class number $h_K$ equals $1$. Analogous to \eqref{eq:correspondanceS2H}, we have a correspondence
        \begin{equation}
        \label{eq:correspondenceNumberField}
            \begin{array}{lll}
                S_2^{++} \left( \mathbb{R} \right)^n & \overset{\sim}{\longrightarrow} & \mathbb{H}^n \times \left( \mathbb{R}_+^{\ast} \right)^n,
            \end{array}
        \end{equation}
        which is compatible with the actions of the extended Hilbert modular group $\widehat{\Gamma}_K$. Using \eqref{eq:correspondenceNumberField}, we can associate, to any point $\tau \in \mathbb{H}^n$, a $2$-dimensional rigid adelic space $E_{\tau}$ with total height $1$. The Roy--Thunder minima
        \begin{equation}
            \begin{array}{lll}
                \Lambda_1 \left( E_{\tau} \right) & = & \inf \left \{ H_{E_{\tau}} \left( x \right) \; \middle \vert \; x \in E_{\tau} \setminus \left \{ 0 \right \} \right \}, \\[0.5em]

                \Lambda_2 \left( E_{\tau} \right) & = & \inf \left \{ \max \left( H_{E_{\tau}} \left( x \right), H_{E_{\tau}} \left( y \right) \right) \; \middle \vert \; \vect_K \left( x,y \right) \, = \, E_{\tau} \right \},
            \end{array}
        \end{equation}
        which can be thought of as analogues of the first and second minima of lattices, fit into a version of Minkowski's second theorem, with the following Hermite-type constant
        \begin{equation}
        \label{eq:intro:HermiteConstant}
            \begin{array}{lll}
                \displaystyle c_{II}^{\Lambda} \left( 2, K \right) & = & \displaystyle \sup\limits_{\dim E = 2} \sqrt{\frac{\Lambda_1 \left( E \right) \Lambda_2 \left( E \right)}{H \left( E \right)}}.
            \end{array}
        \end{equation}
        Similarly to the case of rank $2$ lattices, the correspondence \eqref{eq:correspondenceNumberField} allows us to translate the Roy--Thunder minima of $E_{\tau}$ into invariants linked to the class of $\tau$ in the Hilbert modular variety $\widehat{\Gamma}_K \backslash \mathbb{H}^n$, using the assumption $h_K = 1$. These invariants are linked, respectively, to the distances $\mu_1 \left( \tau \right)^{-1/2}$ to the closest and $\mu_2 \left( \tau \right)^{-1/2}$ to the second closest cusp for the action of $\widehat{\Gamma}_K$. The Minkowski-type theorem which holds for the Roy--Thunder minima then has an incarnation in terms of distances to cusps. It takes the following form.

        \begin{theorem}
            For any $\tau \in \mathbb{H}^n$, we have
            \begin{equation}
                \begin{array}{lllll}
                    \displaystyle \frac{1}{c_{II}^{\Lambda} \left( 2, K \right)^{4n}} & \leqslant & \mu_1 \left( \tau \right) \mu_2 \left( \tau \right) & \leqslant & 1,
                \end{array}
            \end{equation}
            where $c_{II}^{\Lambda} \left( 2, K \right)$ is the constant \eqref{eq:intro:HermiteConstant}.
        \end{theorem}

        Drawing on the equality
        \begin{equation}
            \begin{array}{lll}
                \displaystyle c_{II}^{\Lambda} \left( 2, K \right) & = & \displaystyle \sup\limits_{\dim E = 2} \frac{\Lambda_1 \left( E \right)}{\sqrt{H \left( E \right)}},
            \end{array}
        \end{equation}
        we see that, in the sphere of influence of the cusp $\infty$ in $\mathbb{H}^n$ for the action of $\widehat{\Gamma}_K$, the point $\tau_{max}$ farthest away from $\infty$ is at distance
        \begin{equation}
        \label{eq:distanceTauMax}
            \begin{array}{lll}
                \mu_1 \left( \tau_{\text{max}} \right)^{-1/2} & = & c_{II}^{\Lambda} \left( 2, K \right)^n
            \end{array}
        \end{equation}
        from its closest cusp (\textit{i.e.} $\infty$ in this setting). This is compatible with the observation made in the case of rank $2$ lattices, which corresponds to taking~$K = \mathbb{Q}$. This Minkowski-type theorem has several consequences. The first two are improvements on results in \cite{vanDerGeer:hilbert-modular-surfaces}, under the assumption $h_K \, = \, 1$:
        \begin{enumerate}
            \item corollary \ref{cor:separationCusps} is an effective and uniform version of \cite[Lemma I.2.1]{vanDerGeer:hilbert-modular-surfaces};

            \item corollary \ref{cor:lowerBoundMu1}, is an optimal version of \cite[Lemma I.2.2]{vanDerGeer:hilbert-modular-surfaces}.
        \end{enumerate}
        Note that using the proof of \cite[Lemma I.2.2]{vanDerGeer:hilbert-modular-surfaces} and the optimal nature of corollary~\ref{cor:lowerBoundMu1}, we get
        \begin{equation}
            \begin{array}{lll}
                c_{II}^{\Lambda} \left( 2, K \right) & \leqslant & \sqrt{2} \Delta_K^{1/2n},
            \end{array}
        \end{equation}
        which is the same upper-bound on the Hermite constant $c_{II}^{\Lambda} \left( 2, K \right)$ as the one obtained by Gaudron and Rémond in \cite[Proposition 5.1]{gaudron-remond:corps-siegel} using different methods. The last consequence of the Minkowski-type theorem proved here is related to the problem which motivated this paper: being able to control certain normalized integrals linked to Hilbert modular varieties, which appear in \cite{frey-lefourn-lorenzo:height-estimates}. This result takes the following form.
        \begin{theorem}
            For any $0 \leqslant t < 1$, we have
            \begin{equation}
                \begin{array}{lllll}
                    \multicolumn{5}{l}{\displaystyle \frac{1}{c_{II}^{\Lambda} \left( 2, K \right)^{2nt}} + \frac{t}{1-t} \cdot \frac{1}{c_{II}^{\Lambda} \left( 2, K \right)^{2n}}} \\[2em]
                    
                    \qquad \qquad & \leqslant & \displaystyle \frac{1}{\vol \left( \widehat{\Gamma}_K \backslash \mathbb{H}^n \right)} \; \int_{\widehat{\Gamma}_K \backslash \mathbb{H}^n} \; \mu_1 \left( \tau \right)^t \; \mathrm{d}\mu \left( \tau \right) & \leqslant & \displaystyle \frac{1}{1-t}.
                \end{array}
            \end{equation}
        \end{theorem}
        \noindent In particular, we note that the upper-bound to such integrals can be taken uniformly in $K$, with the requirement that we have $h_K = 1$.

    \subsection{Future work}

        This project has opened up several questions which could be interesting for future work. Below is a (non-exhaustive) list of such problems.

        \begin{enumerate}
            \item The most natural extension would be to drop the requirement $h_K = 1$. The reasoning made in this paper cannot directly be adapted, and one would need to consider more general Hilbert modular groups, as well as the adelic framework presented by van der Geer in \cite[Section I.7]{vanDerGeer:hilbert-modular-surfaces}.

            \item The Hermite constants $c_{II}^{\Lambda} \left( 2, K \right)$ are not explicitly computed in general. It would be interesting to see whether observation \eqref{eq:distanceTauMax} provides a path which could lead to a closed formula. For instance, one may imagine that the point $\tau_{\text{max}}$ has particular properties, and that it could even be an elliptic point, which is the case for $K = \mathbb{Q}$.

            \item Taking into account rigid adelic spaces, or even lattices, of greater rank would also be interesting, and could lead to information on higher rank Hermite constants. The main issue would be to find suitable modular varieties which would fit into a correspondence similar to \eqref{eq:correspondenceNumberField}. Perhaps Siegel modular varieties would be relevant in this context.
        \end{enumerate}

    \subsection{Acknowledgements}

        The question which motivated this work was first presented to me by Samuel Le Fourn, whom I warmly thank for his comments and suggestions. I would also like to thank Gerard Freixas i Montplet, who introduced me to the study of Hilbert modular varieties.

\section{Rigid adelic spaces}

    In this section, we will recall the notion of \textit{rigid adelic spaces}, presented by Gaudron in \cite{gaudron:rigid-adelic-spaces}. Most of this material will be adapted to fit into the framework needed in this paper, which is that of \textit{totally real number fields}.

    \subsection{Algebraic number theory}

        Let us begin by recalling some facts from algebraic number theory. The reader is referred to \cite{neukirch:algebraic-number-theory, ramakrishnan-valenza:fourier-analysis-number-fields} for presentations dedicated to this topic. 

        \subsubsection{Number fields}

            Let $K$ be a number field, \textit{i.e.} a finite field extension of $\mathbb{Q}$, of degree $n$.

            \begin{definition}
                A (multiplicative) valuation\footnote{Also called an absolute value.} of $K$ up to equivalence\footnote{See \cite[Section II.3]{neukirch:algebraic-number-theory} for the equivalence relation.} is called a place. Their set is denoted by $V \left( K \right)$. We further write
                \begin{equation}
                    \begin{array}{lll}
                        V \left( K \right) & = & V_f \left( K \right) \; \sqcup \; V_{\infty} \left( K \right),
                    \end{array}
                \end{equation}
                where $V_f \left( K \right)$ is the set of finite (or non-Archimedean) places of $K$, and $V_{\infty} \left( K \right)$ the set of infinite (or Archimedean) places of $K$. For any place $v \in V \left( K \right)$, we denote by $K_v$ the completion of $K$ with respect to a distance induced by $v$.
            \end{definition}

            \begin{definition}
                We say that $K$ is \textit{totally real} if, for every $v \in V_{\infty} \left( K \right)$, we have an isomorphism of fields $K_v \, \simeq \, \mathbb{R}$, which is the same as having an embedding
                \begin{equation}
                    \begin{array}{ccccc}
                        \sigma & : & K & \hooklongrightarrow & \mathbb{R}.
                    \end{array}
                \end{equation}
                Thus, the fact that $K$ is totally real means we have $n$ different such embeddings, which we denote by $\sigma_1, \ldots, \sigma_n$.
            \end{definition}

            In order to work with valuations on $K$, we need to fix a representative $\left \vert \cdot \right \vert_v$ of each place $v$. This choice is not unique\footnote{Neukirch in \cite[Proposition III.1.2]{neukirch:algebraic-number-theory} and Ramakrishnan--Valenza in \cite[Lemma 5.13]{ramakrishnan-valenza:fourier-analysis-number-fields} use a different convention than Gaudron in \cite{gaudron:rigid-adelic-spaces}, which in turns changes the product formula \eqref{eq:productFormula}.}, and will influence the way certain formulae are written later on.

            \begin{definition}
                Consider two places $v \in V \left( K \right)$ and $w \in V \left( \mathbb{Q} \right)$. We say that $v$ extends $w$, and write $v \vert w$ if the restriction of $v$ to $\mathbb{Q}$ equals $w$.
            \end{definition}

            \begin{proposition}
                Consider a place $v \in V \left( K \right)$, and denote by $w$ its restriction to~$\mathbb{Q}$. There is a unique choice of representative in the equivalence class of $v$ whose restriction to $\mathbb{Q}$ corresponds to either the $p$-adic valuation or the classical absolute value, with the usual normalisations. This choice satisfies
                \begin{equation}
                \label{eq:extensionValuation}
                    \begin{array}{lll}
                        \left \vert \cdot \right \vert_v & = & \left \vert N_{K_v/\mathbb{Q}_w} \left( \cdot \right) \right \vert_w^{1/n_v},
                    \end{array}
                \end{equation}
                with $n_v = \left[ K_v : \mathbb{Q}_w \right]$, and where $N_{K_v/\mathbb{Q}_w} \left( x \right)$ is the determinant of the $\mathbb{Q}_w$-linear map
                \begin{equation}
                    \begin{array}{ccccc}
                        \rho_x & : & K_v & \longmapsto & K_v \\
                        && y & \longmapsto & xy
                    \end{array}.
                \end{equation}
            \end{proposition}

            \begin{proof}
                This is \cite[Lemma 5.13]{ramakrishnan-valenza:fourier-analysis-number-fields}, taking into account the different choice of representative for places of $K$, and can also be found in \cite[Section II.8]{neukirch:algebraic-number-theory}.
            \end{proof}

            \begin{remark}
                Note that, if $v$ is an Archimedean place of $K$, we have $K_v \, = \, \mathbb{Q}_w \, = \, \mathbb{R}$, which gives $n_v \, = \, 1$.
            \end{remark}

            \begin{proposition}[Product formula]
                For any $x \in K^{\times}$, we have $\left \vert x \right \vert_v \, = \, 1$ for all but finitely many places $v \in V \left( K \right)$, as well as
                \begin{equation}
                \label{eq:productFormula}
                    \begin{array}{lll}
                        \prod\limits_{v \in V \left( K \right)} \left \vert x \right \vert_v^{n_v} & = & 1.
                    \end{array}
                \end{equation}
            \end{proposition}

            \begin{proof}
                Consider $x \in K^{\times}$. We have
                \begin{equation}
                    \begin{array}{lllll}
                        \prod\limits_{v \in V \left( K \right)} \left \vert x \right \vert_v^{n_v} & = & \prod\limits_{w \in V \left( \mathbb{Q} \right)} \prod\limits_{v \vert w} \left \vert x \right \vert_v^{n_v} & = & \prod\limits_{w \in V \left( \mathbb{Q} \right)} \prod\limits_{v \vert w} \left \vert N_{K_v/\mathbb{Q}_w} \left( x \right) \right \vert_w \\[2em]

                        &&& = & \prod\limits_{w \in V \left( \mathbb{Q} \right)} \left \vert N_{K/\mathbb{Q}} \left( x \right) \right \vert_w \\[2em]

                        &&& = & 1,
                    \end{array}
                \end{equation}
                using the classical product formula on $\mathbb{Q}$. Note that we have also used the canonical isomorphism (see \cite[Proposition II.8.3]{neukirch:algebraic-number-theory})
                \begin{equation}
                    \begin{array}{lll}
                        K \otimes \mathbb{Q}_w & \simeq & \prod\limits_{v \vert w} K_v
                    \end{array}
                \end{equation}
                for any $w \in V \left( \mathbb{Q} \right)$ to get the equality
                \begin{equation}
                    \begin{array}{lll}
                        \prod\limits_{v \vert w} \left \vert N_{K_v/\mathbb{Q}_w} \left( x \right) \right \vert_w & = & \left \vert N_{K/\mathbb{Q}} \left( x \right) \right \vert_w.
                    \end{array}
                \end{equation}
            \end{proof}

            \begin{definition}
                The ring of ad\`{e}les of $K$ is defined as
                \begin{equation}
                    \begin{array}{lll}
                        \mathbb{A}_K & = & \left \{ \left( x_v \right) \in \prod\limits_{v \in V \left( K \right)} K_v \; \middle \vert \; \text{for almost all } v \in V_f \left( K \right), \; x_v \in \mathcal{O}_{K_v} \right \}.
                    \end{array}
                \end{equation}
            \end{definition}

            \begin{definition}
                The \textit{adelic norm} $\left \vert \cdot \right \vert_{\mathbb{A}_K}$ is defined by
                \begin{equation}
                    \begin{array}{lll}
                        \left \vert x \right \vert_{\mathbb{A}_K} & = & \prod\limits_{v \in V \left( K \right)} \left \vert x_v \right \vert_v^{n_v}
                    \end{array}
                \end{equation}
                for any $x \, = \, \left( x_v \right)_{v \in V \left( K \right)} \in \mathbb{A}_K$.
            \end{definition}

            \begin{proposition}
            \label{prop:normUnit}
                For any $\varepsilon \in \mathcal{O}_K^{\times}$, we have
                \begin{equation}
                    \begin{array}{lllll}
                        N \left( \varepsilon \right) & = & \prod\limits_{i=1}^n \sigma_i \left( \varepsilon \right) & \in & \left \{ \pm 1 \right\}.
                    \end{array}
                \end{equation}
            \end{proposition}

            \begin{proof}
                This is a consequence of the product formula \eqref{eq:productFormula}, and of the fact that a unit $\varepsilon$ satisfies $\left \vert \varepsilon \right \vert_v \, = \, 1$ for any $v \in V_f \left( K \right)$.
            \end{proof}

        \subsubsection{Some invariants of number fields}

            In this paragraph, we will recall the definition of some of the invariants which will appear in later sections.

            \begin{definition}
                Consider a $\mathbb{Z}$-basis $e_1, \ldots, e_n$ of $\mathcal{O}_K$. The \textit{discriminant} of $K$, denoted by $\Delta_K$ is defined as
                \begin{equation}
                    \begin{array}{lll}
                        \Delta_K & = & \det \left( \left( \sigma_i e_j \right)_{i,j = 1, \ldots, n} \right)^2,
                    \end{array}
                \end{equation}
                where we recall that $\sigma_1, \ldots, \sigma_n$ are the real embeddings corresponding to the Archimedean places of $K$. It does not depend on the choice of $\mathbb{Z}$-basis of $\mathcal{O}_K$.
            \end{definition}

            \begin{proposition}
            \label{prop:actionOKRn}
                The map
                \begin{equation}
                    \begin{array}{ccccc}
                        \sigma & : & \mathcal{O}_K & \longrightarrow & \mathbb{R}^n \\[0.5em]
                        && \alpha & \longmapsto & \left( \sigma_1 \left( \alpha \right), \ldots, \sigma_n \left( \alpha \right) \right)
                    \end{array}
                \end{equation}
                embeds $\mathcal{O}_K$ as a full-rank lattice in $\mathbb{R}^n$, and we have
                \begin{equation}
                    \begin{array}{lll}
                        \vol \left( \mathcal{O}_K \backslash \mathbb{R}^n \right) & = & \sqrt{\Delta_K}.
                    \end{array}
                \end{equation}
            \end{proposition}

            \begin{proof}
                This is a direct consequence of the definition of the discriminant.
            \end{proof}

            \begin{definition}
                Denote by $\mathcal{O}_K^{\times, +}$ the group of totally positive units of $K$, \textit{i.e.} of elements $\eta \in \mathcal{O}_K^{\times}$ such that we have $\sigma_i \left( \eta \right) > 0$ for all integers $i \in \llbracket 1,n \rrbracket$.
            \end{definition}
            
            \begin{definition}
                We define the group $\mathcal{O}_K^{\times,2}$ as
                \begin{equation}
                    \begin{array}{lll}
                        \mathcal{O}_K^{\times,2} & = & \left \{ x^2 \; \middle \vert \; x \in \mathcal{O}_K^{\times} \right \}.
                    \end{array}
                \end{equation}
            \end{definition}

            \begin{remark}
                We have the inclusions
                \begin{equation}
                \label{eq:inclusionsOK}
                    \begin{array}{ccccc}
                        \mathcal{O}_K^{\times, 2} & \subseteq & \mathcal{O}_K^{\times, +} & \subseteq & \mathcal{O}_K^{\times}.
                    \end{array}
                \end{equation}
            \end{remark}

            \begin{theorem}[Dirichlet's unit theorem]
            \label{thm:dirichetUnit}
                There exists a group isomorphism
                \begin{equation}
                    \label{eq:dirichletUnit}
                    \begin{array}{lll}
                        \mathcal{O}_K^{\times} & \simeq & \left( \mathbb{Z}/2\mathbb{Z} \right) \times \mathbb{Z}^{n-1}.
                    \end{array}
                \end{equation}
            \end{theorem}

            \begin{proof}
                This is \cite[Theorem I.7.4]{neukirch:algebraic-number-theory}, keeping in mind that $K$ is totally real, which means that the roots of unity can only be $\pm 1$.
            \end{proof}

            \begin{corollary}
                We have
                \begin{equation}
                    \begin{array}{lll}
                        [ \mathcal{O}_K^{\times} : \mathcal{O}_K^{\times,2} ] & = & 2^n.
                    \end{array}
                \end{equation}
            \end{corollary}

            \begin{proof}
                This is a direct consequence of theorem \ref{thm:dirichetUnit}, noting that, under the isomorphism \eqref{eq:dirichletUnit}, we have
                \begin{equation}
                    \begin{array}{lll}
                        \mathcal{O}_K^{\times, 2} & \simeq & \left( 2 \mathbb{Z} \right)^{n-1}.
                    \end{array}
                \end{equation}
            \end{proof}

            Consider a set $\varepsilon_1, \ldots, \varepsilon_{n-1}$ of generators of $\mathcal{O}_K^{\times} / \left \{ \pm 1 \right \}$.

            \begin{definition}
                The \textit{regulator} of $K$, denoted by $R_K$, is defined as the absolute value of any rank $n-1$ minor of the matrix
                \begin{equation}
                \label{eq:matrixRegulator}
                    \begin{array}{c}
                        \left( \begin{array}{ccc} \log \sigma_1 \left( \varepsilon_1 \right) & \ldots & \log \sigma_1 \left( \varepsilon_{n-1} \right) \\ \vdots & & \vdots \\ \log \sigma_n \left( \varepsilon_1 \right) & \ldots & \log \sigma_n \left( \varepsilon_{n-1} \right) \end{array} \right).
                    \end{array}
                \end{equation}
                The result does not depend on the choice of generators $\varepsilon_1, \ldots, \varepsilon_{n-1}$.
            \end{definition}

            \begin{remark}
                The fact that any rank $n-1$ minor of the matrix \eqref{eq:matrixRegulator} can be taken to define the regulator comes from proposition \ref{prop:normUnit}, which implies that, for any $\varepsilon \in \mathcal{O}_K^{\times, +}$, we have
                \begin{equation}
                \label{eq:normOfTotPosUnit}
                    \begin{array}{lllll}
                        N \left( \varepsilon \right) & = & \prod\limits_{i=1}^n \sigma_i \left( \varepsilon \right) & = & 1.
                    \end{array}
                \end{equation}
            \end{remark}

            \begin{proposition}
            \label{prop:OKTimesPlusLatticeRNMinus1}
                Using $\sigma_1, \ldots, \sigma_{n-1}$, the map
                \begin{equation}
                \label{eq:mapLambdaNMinus1}
                    \begin{array}{ccccc}
                        \lambda & : & \left( \mathbb{R}_+^{\ast} \right)^{n-1} & \longrightarrow & \mathbb{R}^{n-1} \\[0.5em]
                        && \left( y_1, \ldots, y_{n-1} \right) & \longmapsto & \left( \log y_1, \ldots, \log y_{n-1} \right)
                    \end{array}
                \end{equation}
                embeds $\mathcal{O}_K^{\times,+}$ as a full-rank lattice in $\mathbb{R}^{n-1}$, whose volume is given by
                \begin{equation}
                \label{eq:volumeActionOKTimesPlusRNMinus1}
                    \begin{array}{lll}
                        \vol \left( \mathcal{O}_K^{\times,+} \backslash \mathbb{R}^{n-1} \right) & = & \displaystyle \frac{2^{n-1}}{[ \mathcal{O}_K^{\times,+} : \mathcal{O}_K^{\times,2}]} \cdot R_K.
                    \end{array}
                \end{equation}
            \end{proposition}

            \begin{proof}
                The map \eqref{eq:mapLambdaNMinus1} also embeds $\mathcal{O}_K^{\times} / \left \{ \pm 1 \right \}$ as a full-rank lattice in $\mathbb{R}^{n-1}$, and, by definition of the regulator, we have
                \begin{equation}
                    \begin{array}{lll}
                        \vol \left( \left( \mathcal{O}_K^{\times} / \left \{ \pm 1 \right \} \right) \backslash \mathbb{R}^{n-1} \right) & = & R_K.
                    \end{array}
                \end{equation}
                The remaining factor in \eqref{eq:volumeActionOKTimesPlusRNMinus1} comes from the equality
                \begin{equation}
                    \begin{array}{ccccc}
                        [ \mathcal{O}_K^{\times} / \left \{ \pm 1 \right \} : \mathcal{O}_K^{\times,+}] & = & \displaystyle \frac{[ \mathcal{O}_K^{\times} / \left \{ \pm 1 \right \} : \mathcal{O}_K^{\times,2}]}{[ \mathcal{O}_K^{\times,+} : \mathcal{O}_K^{\times,2}]} & = & \displaystyle \frac{2^{n-1}}{[ \mathcal{O}_K^{\times,+} : \mathcal{O}_K^{\times,2}]}.
                    \end{array}
                \end{equation}
            \end{proof}

            Before moving on to the next invariant, note that
            \begin{equation}
            \label{eq:actionOKTimesRnPlus}
                \begin{array}{lll}
                    \varepsilon \cdot \left( y_1, \ldots, y_n \right) & = & \left( \sigma_1 \left( \varepsilon \right) y_1, \ldots, \sigma_n \left( \varepsilon \right) \right)
                \end{array}
            \end{equation}
            induces an action of $\mathcal{O}_K^{\times,+}$ on $\left( \mathbb{R}_+^{\ast} \right)^n$, which, using \eqref{eq:normOfTotPosUnit}, preserves
            \begin{equation}
                \begin{array}{lll}
                    B & = & \left \{ \left( y_1, \ldots, y_n \right) \in \left( \mathbb{R}_+^{\ast} \right)^n \; \middle \vert \; y_1 \ldots y_n \, = \, 1 \right \}.
                \end{array}
            \end{equation}
            We will need to compute the volume of a fundamental domain for this action in the next section, where the measure used on $\left( \mathbb{R}_+^{\ast} \right)^n$ will appear naturally.

            We will now define the \textit{class number} of $K$, omitting many of the details, for which the reader is referred to \cite{neukirch:algebraic-number-theory}.

            \begin{proposition-definition}
                Denote by $J_K$ the group of fractional ideals of $K$, and by $P_K \subseteq J_K$ the subgroup of principal fractional ideals. The quotient group
                \begin{equation}
                    \begin{array}{lll}
                        Cl_K & = & J_K / P_K
                    \end{array}
                \end{equation}
                is finite, and its order $h_K$ is called the class number of $K$.
            \end{proposition-definition}

            \begin{proof}
                The finiteness of $Cl_K$ is proved in \cite[Theorem I.6.3]{neukirch:algebraic-number-theory}.
            \end{proof}

            \begin{definition}
                The \textit{Dedekind zeta function} $\zeta_K$ of $K$ is defined as
                \begin{equation}
                \label{eq:dedekindZeta}
                    \begin{array}{ccccc}
                        \zeta_K & : & s & \longmapsto & \sum\limits_{\mathfrak{a}} N \left( \mathfrak{a} \right)^{-s},
                    \end{array}
                \end{equation}
                for any $s \in \mathbb{C}$ with $\Real s > 1$, and where the sum in \eqref{eq:dedekindZeta} ranges over ideals $\mathfrak{a} \subseteq \mathcal{O}_K$, with $N \left( \mathfrak{a} \right) = \# \left( \mathcal{O}_K / \mathfrak{a} \right)$.
            \end{definition}

            \begin{proposition}
                The Dedekind zeta function admits a meromorphic continuation to $\mathbb{C}$ with only one pole at $s=1$, which is of order $1$.
            \end{proposition}

            \begin{proof}
                This is done in \cite[Corollary VII.5.11]{neukirch:algebraic-number-theory}.
            \end{proof}

    \subsection{Rigid adelic spaces}

        In this section, we will present the work of Gaudron from \cite{gaudron:rigid-adelic-spaces}, in the context of totally real number fields.

        \begin{definition}
            An \textit{adelic space} over $K$ is the datum of a finite-dimensional vector space $E$ over $K$, together with a norm $\left \Vert \cdot \right \Vert_{E, v}$ on each $E \otimes_K K_v$, with $v \in V \left( K \right)$.
        \end{definition}

        From now on, we will restrict ourselves to the dimension $2$ case.

        \begin{definition}
            The \textit{standard adelic space} $E_0$ of dimension $2$ is $K^2$ endowed with the norms
            \begin{equation}
                \begin{array}{lll}
                    \left \Vert \left( x_1, x_2 \right) \right \Vert_v & = & \left \{ \begin{array}{ll} \max \left( \left \vert x_1 \right \vert_v, \left \vert x_2 \right \vert_v \right) & \text{if } v \in V_f \left( K \right) \\[1em] \sqrt{\left \vert \sigma_i \left( x_1 \right) \right \vert^2 + \left \vert \sigma_i \left( x_2 \right) \right \vert^2} & \text{if } v \in V_{\infty} \left( K \right) \end{array} \right.
                \end{array},
            \end{equation}
            keeping in mind that an infinite place is associated with exactly one of the real embeddings $\sigma_1, \ldots, \sigma_n$.
        \end{definition}

        \begin{definition}
            Define the group $GL_2 \left( \mathbb{A}_K \right)$ by
            \begin{equation}
                \begin{array}{lllll}
                     GL_2 \left( \mathbb{A}_K \right) & = & \left \{ \left( A_v \right)_{v \in V \left( K \right)} \in \prod\limits_{v \in V \left( K \right)} GL_2 \left( K_v \right) \; \right \vert \\[2em]
                     
                     && \qquad \quad \left. \text{for almost all } v \in V_f \left( K \right), \; A_v \, \in \, GL_2 \left( \mathcal{O}_{K_v} \right) \right \}.
                \end{array}
            \end{equation}
        \end{definition}

        \begin{definition}
        \label{def:rigidAdelicSpaceDim2}
            A \textit{rigid adelic space} $E$ of dimension $2$ is the datum of an adelic space $E$ of dimension $2$ such that, for every $v \in V \left( K \right)$, there exists a linear isomorphism
            \begin{equation}
                \begin{array}{ccccc}
                    \varphi_v & : & E \otimes_K K_v & \longrightarrow & K_v^2,
                \end{array}
            \end{equation}
            and an element $A = \left( A_v \right)_{v \in V \left( K \right)} \in GL_2 \left( \mathbb{A}_K \right)$, satisfying
            \begin{equation}
                \begin{array}{lll}
                    \left \Vert \left( x_1, x_2 \right) \right \Vert_{E,v} & = & \left \Vert A_v \phi_v \left( x_1, x_2 \right) \right \Vert_v.
                \end{array}
            \end{equation}
        \end{definition}

    \subsection{Heights and Roy--Thunder minima}

        To any rigid adelic space, one can attach several invariants, which the reader can find in \cite{gaudron:rigid-adelic-spaces, gaudron-remond:corps-siegel}. Let us present the ones which will be useful in this paper, in the setting we are concerned with. Below, we will consider a dimension $2$ rigid adelic space $E$ associated with $A \in GL_2 \left( \mathbb{A}_K \right)$.

        \begin{definition}
        \label{def:HeightTotal}
            The \textit{height} of $E$ is defined as
            \begin{equation}
                \begin{array}[t]{lllll}
                    H \left( E \right) & = & \left \vert \det A \right \vert_{\mathbb{A}_K}^{1/n} & = & \prod\limits_{v \in V \left( K \right)} \left \vert \det A_v \right \vert_v^{n_v/n}.
                \end{array}
            \end{equation}
        \end{definition}

        \begin{definition}
        \label{def:HeightPoint}
            Let $x \in E$. The \textit{height of $x$} is defined as
            \begin{equation}
                \begin{array}[t]{lll}
                    H_E \left( x \right) & = & \prod\limits_{v \in V \left( K \right)} \left \Vert x \right \Vert_{E,v}^{n_v/n}.
                \end{array}
            \end{equation}
        \end{definition}

        \begin{proposition}
            The height function $H_E$ is projective, insofar as we have
            \begin{equation}
                \begin{array}{lll}
                    H_E \left( \lambda x \right) & = & H_E \left( x \right)
                \end{array}
            \end{equation}
            for any $\lambda \in K^{\ast}$ and any $x \in E$.
        \end{proposition}

        \begin{proof}
            This is a direct consequence of the product formula \eqref{eq:productFormula}.
        \end{proof}

        \begin{definition}
            The \textit{Roy--Thunder minima} of $E$ are defined as
            \begin{equation}
                \begin{array}{lll}
                    \Lambda_1 \left( E \right) & = & \inf \left \{ H_E \left( x \right) \; \middle \vert \; x \in E \setminus \left \{ 0 \right \} \right \}, \\[0.5em]

                    \Lambda_2 \left( E \right) & = & \inf \left \{ \max \left( H_E \left( x \right), H_E \left( y \right) \right) \; \middle \vert \; \vect_K \left( x,y \right) \, = \, E \right \}.
                \end{array}
            \end{equation}
        \end{definition}

        \begin{definition}
        \label{def:HermiteConstant}
            The Hermite constant $c_{II}^{\Lambda} \left( 2, K \right)$ is defined as
            \begin{equation}
                \begin{array}{lll}
                    \displaystyle c_{II}^{\Lambda} \left( 2, K \right) & = & \displaystyle \sup\limits_{\dim E = 2} \sqrt{\frac{\Lambda_1 \left( E \right) \Lambda_2 \left( E \right)}{H \left( E \right)}},
                \end{array}
            \end{equation}
            where the supremum ranges over all rigid adelic spaces of dimension $2$.
        \end{definition}

        \begin{remark}
            These Hermite constants are always finite, since $K$ is a number field (see \cite[after definition 25]{gaudron:rigid-adelic-spaces}).
        \end{remark}

        \begin{proposition}
        \label{prop:MinkowskiRoyThunder}
            For any $2$-dimensional rigid adelic space $E$, we have
            \begin{equation}
                \begin{array}{lllll}
                    H \left( E \right) & \leqslant & \Lambda_1 \left( E \right) \Lambda_2 \left( E \right) & \leqslant & H \left( E \right) c_{II}^{\Lambda} \left( 2, K \right)^2.
                \end{array}
            \end{equation}
        \end{proposition}

        \begin{proof}
            The right-hand side inequality is a direct consequence of the definition of the Hermite constant $c_{II}^{\Lambda} \left( 2, K \right)$. The left-hand side one is a version of Hadamard's inequality, as noted in \cite[Lemme 3.7]{gaudron-remond:corps-siegel}.
        \end{proof}

\section{\texorpdfstring{Action of Hilbert modular groups on $\mathbb{H}^n$}{Action of Hilbert modular groups on Hn}}

    This section will be devoted to a quick overview of the geometric objects used in this paper, which are \textit{Hilbert modular varieties}. The reader will find ample information about these in \cite{hirzebruch:hilbert-modular, vanDerGeer:hilbert-modular-surfaces}, even though these references mainly deal with Hilbert modular surfaces. Once again, we consider a totally real number field $K$ of degree $n$, and use the notations introduced in the previous section.

    \subsection{Hilbert modular groups}

        The groups we will be concerned with can be seen as generalizations of the group of unimodular matrices $PSL_2 \left( \mathbb{Z} \right)$ to totally real number fields, and act on the product of $n$ copies of the Poincaré upper-half plane $\mathbb{H} \, = \, \left\{ x+iy \in \mathbb{C} \; \middle \vert \; y > 0 \right \}$.

        \subsubsection{Action by homographic transformations}

            Consider the group
            \begin{equation}
                \begin{array}{lll}
                    PGL_2^+ \left( \mathbb{R} \right) & = & \left \{ \left( \begin{array}{cc} a & b \\[0.4em] c & d \end{array} \right) \in PGL_2 \left( \mathbb{R} \right) \; \middle \vert \; ad-bc>0 \right \},
                \end{array}
            \end{equation}
            which acts on the Riemann sphere $\mathbb{P}^1 \left( \mathbb{C} \right)$ in the usual way, via
            \begin{equation}
            \label{eq:actionHomographic}
                \begin{array}{lll}
                    \left( \begin{array}{cc} a & b \\[0.4em] c & d \end{array} \right) \cdot \tau & = & \displaystyle \frac{a \tau + b}{c \tau + d}.
                \end{array}
            \end{equation}
            Writing
            \begin{equation}
                \begin{array}{lll}
                    \mathbb{P}^1 \left( \mathbb{C} \right) & = & \mathbb{H} \; \sqcup \; \mathbb{P}^1 \left( \mathbb{R} \right) \; \sqcup \; \overline{\mathbb{H}},
                \end{array}
            \end{equation}
            where $\overline{\mathbb{H}}$ is the complex conjugate of the Poincaré upper-half plane (\textit{i.e.} the lower-half plane), we note that the action \eqref{eq:actionHomographic} preserves $\mathbb{H}$, $\mathbb{P}^1 \left( \mathbb{R} \right)$, and $\overline{\mathbb{H}}$.
            \begin{proposition}
                An element $\gamma \in PGL_2^+ \left( \mathbb{R} \right)$ different from the identity has either:
                \begin{enumerate}
                    \item a single fixed point in $\mathbb{P}^1 \left( \mathbb{R} \right)$, in which case $\gamma$ is said to be parabolic;
                    \item two distinct fixed points in $\mathbb{P}^1 \left( \mathbb{R} \right)$, in which case $\gamma$ is said to be hyperbolic;
                    \item two complex conjugated fixed points in $\mathbb{H}$ and $\overline{\mathbb{H}}$, in which case $\gamma$ is said to be elliptic.
                \end{enumerate}
            \end{proposition}

            \begin{proof}
                This is a direct computation.
            \end{proof}

            By considering $n$ copies of this, we get an action of $PGL_2^+ \left( \mathbb{R} \right)^n$ on $\mathbb{H}^n$, $\mathbb{P}^1 \left( \mathbb{R} \right)^n$, and $\overline{\mathbb{H}}^n$.

            \begin{definition}
                An element $\gamma \, = \, \left( \gamma_1, \ldots, \gamma_n \right) \in PGL_2^+ \left( \mathbb{R} \right)^n$ is said to be:
                \begin{enumerate}
                    \item \textit{parabolic} if every $\gamma_i$ is parabolic;
                    \item \textit{hyperbolic} if every $\gamma_i$ is hyperbolic;
                    \item \textit{elliptic} if every $\gamma_i$ is elliptic;
                    \item \textit{mixed} otherwise.
                \end{enumerate}
            \end{definition}

        \subsubsection{The Hilbert modular groups}

            There are  be two subgroups of $PGL_2 \left( \mathbb{R} \right)^n$ which will be of particular importance in this paper.

            \begin{definition}
                The \textit{Hilbert modular group} $\Gamma_K$ is defined as
                \begin{equation}
                    \begin{array}{lllll}
                        \Gamma_K & = & PSL_2 \left( \mathcal{O}_K \right) & = & SL_2 \left( \mathcal{O}_K \right) / \left \{ \pm 1 \right \},
                    \end{array}
                \end{equation}
                and the \textit{extended Hilbert modular group} $\widehat{\Gamma}_K$ as
                \begin{equation}
                    \begin{array}{lll}
                        \widehat{\Gamma}_K & = & \left \{ M \in GL_2 \left( \mathcal{O}_K \right) \; \middle \vert \; \det M \in \mathcal{O}_K^{\times, +} \right \} / \left \{ \left( \begin{array}{cc} \varepsilon & 0 \\ 0 & \varepsilon \end{array} \right) \; \middle \vert \; \varepsilon \in \mathcal{O}_K^{\times} \right \}.
                    \end{array}
                \end{equation}
            \end{definition}

            \begin{proposition}
                We have an exact sequence of groups
                \begin{equation}
                    \begin{array}{ccccccccc}
                        1 & \longrightarrow & \Gamma_K & \longrightarrow & \widehat{\Gamma}_K & \longrightarrow & \mathcal{O}_K^{\times,+}/\mathcal{O}_K^{\times, 2} & \longrightarrow & 1.
                    \end{array}
                \end{equation}
                As a consequence, we have
                \begin{equation}
                \label{eq:indexHilbertModularGroups}
                    \begin{array}{lll}
                        [ \widehat{\Gamma}_K : \Gamma_K ] & = & [ \mathcal{O}_K^{\times, +} : \mathcal{O}_K^{\times, 2} ].
                    \end{array}
                \end{equation}
            \end{proposition}

            \begin{proof}
                The injective group homomorphism
                \begin{equation}
                    \begin{array}{lll}
                        \Gamma_K & \hooklongrightarrow & \widehat{\Gamma}_K
                    \end{array}
                \end{equation}
                is the natural inclusion which comes from the definition of the Hilbert modular groups, while the surjective group homomorphism
                \begin{equation}
                    \begin{array}{lll}
                        \widehat{\Gamma}_K & \twoheadlongrightarrow & \mathcal{O}_K^{\times,+}/\mathcal{O}_K^{\times, 2}
                    \end{array}
                \end{equation}
                is given by the determinant. The rest of the proof is a direct computation.
            \end{proof}

            Using the real embeddings $\sigma_1, \ldots, \sigma_n$, we now have
            \begin{equation}
                \begin{array}{lllll}
                    \Gamma_K & \subseteq & \widehat{\Gamma}_K & \subset & PGL_2^+ \left( \mathbb{R} \right)^n,
                \end{array}
            \end{equation}
            so the Hilbert modular groups act on $\mathbb{H}^n$ and on $\mathbb{P}^1 \left( \mathbb{R} \right)^n$.

            \begin{remark}
                The quotients $\Gamma_K \backslash \mathbb{H}^n$ and $\widehat{\Gamma}_K \backslash \mathbb{H}^n$ are examples of \textit{Hilbert modular varieties}, which are noncompact complex analytic spaces.
            \end{remark}

        \subsubsection{The Poincaré metric}

            In order to measure distances and volumes in a meaningful way on $\mathbb{H}^n$, we need to consider a metric. This is done by generalizing the Poincaré metric on $\mathbb{H}$ to the higher-dimensional setting.

            \begin{definition}
                The \textit{Poincaré metric} on $\mathbb{H}^n$ is defined by
                \begin{equation}
                    \begin{array}{lll}
                        \mathrm{d}\mu \left( \tau \right) & = & \displaystyle \frac{\mathrm{d}x_1^2 + \mathrm{d}y_1^2}{y_1^2} \cdot \ldots \cdot \frac{\mathrm{d}x_n^2 + \mathrm{d}y_n^2}{y_n^2}
                    \end{array}
                \end{equation}
                where $\tau \in \mathbb{H}^n$ is written as
                \begin{equation}
                    \begin{array}{lllll}
                        \tau & = & \left( \tau_1, \ldots, \tau_n \right) & = & \left( x_1 + i y_1, \ldots, x_n + i y_n \right).
                    \end{array}
                \end{equation}
            \end{definition}

            \begin{proposition}
                The Poincaré metric is invariant by the action of $\widehat{\Gamma}_K$.
            \end{proposition}

            \begin{proof}
                The computation is entirely similar to the $1$-dimensional case.
            \end{proof}

            \begin{proposition}
                The volume of the Hilbert modular variety $\Gamma_K \backslash \mathbb{H}^n$, \textit{i.e.} of any fundamental domain for the action of $\Gamma_K$ on $\mathbb{H}^n$, with respect to the Poincaré metric is given by
                \begin{equation}
                \label{eq:volumeHMV}
                    \begin{array}{lll}
                        \vol \left( \Gamma_K \backslash \mathbb{H}^n \right) & = & \displaystyle \frac{2}{\pi^n} \cdot \Delta_K^{3/2} \cdot \zeta_K \left( 2 \right).
                    \end{array}
                \end{equation}
            \end{proposition}

            \begin{proof}
                This computation is done by Siegel\footnote{Note there is a gap in \cite{siegel:volume-fundamental-domain}, fixed by Siegel in \cite{siegel:bestimmung-volumens}.} in \cite{siegel:volume-fundamental-domain, siegel:bestimmung-volumens}.
            \end{proof}

            \begin{remark}
                The volume of the Hilbert modular variety $\widehat{\Gamma}_K \backslash \mathbb{H}^n$ is obtained from formula \eqref{eq:volumeHMV} by dividing by $[ \widehat{\Gamma}_K : \Gamma_K ]$, whose expression is given in \eqref{eq:indexHilbertModularGroups}.
            \end{remark}

    \subsection{Cusps}

        The action of the Hilbert modular groups on $\mathbb{P}^1 \left( \mathbb{R} \right)^n$ leads to the definition of specials points, called cusps, which are typically used to compactify Hilbert modular varieties.

        \subsubsection{The set of cusps}

            \begin{definition}
                The set of \textit{cusps} for the action of $\widehat{\Gamma}_K$ is defined as
                \begin{equation}
                    \begin{array}{lll}
                        \left \{ c \in \mathbb{P}^1 \left( \mathbb{R} \right)^n \; \middle \vert \; \exists \, \gamma \in \widehat{\Gamma}_K \text{ parabolic }, \; \gamma \cdot c \, = \, c \right \}.
                    \end{array}
                \end{equation}
            \end{definition}

            Using the real embeddings $\sigma_1, \ldots, \sigma_n$, the projective line $\mathbb{P}^1 \left( K \right)$ is embedded in~$\mathbb{P}^1 \left( \mathbb{R} \right)^n$.

            \begin{proposition}
                The set of cusps of $\widehat{\Gamma}_K$ equals $\mathbb{P}^1 \left( K \right)$, and this set is preserved by the action of the Hilbert modular groups.
            \end{proposition}

            \begin{proof}
                This is a direct computation.
            \end{proof}

            \begin{definition}
                For any cusp $c \in \mathbb{P}^1 \left( K \right)$, we define
                \begin{equation}
                    \begin{array}{lll}
                        \widehat{\Gamma}_{K,c} & = & \left \{ \gamma \in \widehat{\Gamma}_K \; \middle \vert \; \gamma \cdot c \, = \, c \right \},
                    \end{array}
                \end{equation}
                \textit{i.e.} the stabilizer of $c$, which is a non-trivial subgroup of $\widehat{\Gamma}_K$.
            \end{definition}

            \begin{proposition}
                The map
                \begin{equation}
                    \begin{array}{ccc}
                        \widehat{\Gamma}_K \backslash \mathbb{P}^1 \left( K \right) & \longrightarrow & Cl_K \\[0.5em]
                        \left[ \alpha : \beta \right] & \longmapsto & \alpha \mathcal{O}_K + \beta \mathcal{O}_K
                    \end{array}
                \end{equation}
                is well-defined and bijective. Consequently, there are exactly $h_K$ cusps up to $\widehat{\Gamma}_K$-equivalence.
            \end{proposition}

            \begin{proof}
                This is \cite[Proposition I.1.1]{vanDerGeer:hilbert-modular-surfaces}.
            \end{proof}

            \begin{remark}
                There are also $h_K$ cusps up to $\Gamma_K$-equivalence.
            \end{remark}

            \begin{definition}
                The cusp $\infty$ is defined as $\left[ 1 : 0 \right] \in \mathbb{P}^1 \left( K \right)$.
            \end{definition}

            \begin{proposition}
                The stabilizer $\widehat{\Gamma}_{K, \infty}$ of the cusp $\infty$ is given by
                \begin{equation}
                    \begin{array}{lll}
                        \widehat{\Gamma}_{K, \infty} & = & \left \{ \left( \begin{array}{cc} \varepsilon & \mu \\[0.4em] 0 & 1 \end{array} \right) \; \middle \vert \; \varepsilon \in \mathcal{O}_K^{\times, +}, \; \mu \in \mathcal{O}_K \right \}.
                    \end{array}
                \end{equation}
            \end{proposition}

            \begin{proof}
                This is done in \cite[Section I.4]{vanDerGeer:hilbert-modular-surfaces}.
            \end{proof}

            In the following, to lighten notations, for any $\tau \, = \, \left( \tau_1, \ldots, \tau_n \right) \in \mathbb{H}^n$, we set
            \begin{equation}
                \begin{array}{ccccc}
                    \Real \tau & = & \left( \Real \tau_1, \ldots, \Real \tau_n \right) & \in & \mathbb{R}^n, \\[0.5em]
                    \Imag \tau & = & \left( \Imag \tau_1, \ldots, \Imag \tau_n \right) & \in & \left( \mathbb{R}_+^{\ast} \right)^n.
                \end{array}
            \end{equation}

            \begin{proposition}
            \label{prop:fundamentalDomainStabInfty}
                A fundamental domain for the action of $\widehat{\Gamma}_{K, \infty}$ on $\mathbb{H}^n$ is given by the set
                \begin{equation}
                \label{eq:fundamentalDomainStabInfty}
                    \begin{array}{lllll}
                        E & = & \left \{ \tau \in \mathbb{H}^n \; \middle \vert \; \Real \tau \in T, \; \Imag \tau \in F \right \} & = & T \times F,
                    \end{array}
                \end{equation}
                where $T$ is a fundamental domain for the action of $\mathcal{O}_K$ on $\mathbb{R}^n$, as considered in proposition \ref{prop:actionOKRn}, and $F$ is a fundamental domain for the action of $\mathcal{O}_K^{\times, +}$ on $\left( \mathbb{R}_+^{\ast} \right)^n$ considered in \eqref{eq:actionOKTimesRnPlus}.
            \end{proposition}

            \begin{proof}
                For any $\varepsilon \in \mathcal{O}_K^{\times, +}$ and $\mu \in \mathcal{O}_K$, recall that we have
                \begin{equation}
                    \begin{array}{lll}
                        \left( \begin{array}{cc} \varepsilon & \mu \\[0.4em] 0 & 1 \end{array} \right) \cdot \tau & = & \varepsilon \tau + \mu.
                    \end{array}
                \end{equation}
                Starting with a point $\tau \in \mathbb{H}^n$, we begin by setting $\varepsilon \in \mathcal{O}_K^{\times,+}$ so that we have
                \begin{equation}
                    \begin{array}{lll}
                        \Imag \left( \varepsilon \tau \right) & \in & F,
                    \end{array}
                \end{equation}
                and we then fix $\mu \in \mathcal{O}_K$ so as to have
                \begin{equation}
                    \begin{array}{lll}
                        \Real \left( \varepsilon \tau + \mu \right) & \in & T.
                    \end{array}
                \end{equation}
                Noting that we have $\Imag \left( \varepsilon \tau + \mu \right) \, = \, \Imag \left( \varepsilon \tau \right)$, any point in $\mathbb{H}^n$ can be brought to $E$ by applying a matrix in $\widehat{\Gamma}_{K, \infty}$. A similar reasoning proves that $\tau \in E$ cannot be sent by $\widehat{\Gamma}_{K,\infty}$ to a different point of $E$.
            \end{proof}

            To conclude this paragraph, let us find a convenient example of fundamental domain $F$ which we can plug into \eqref{eq:fundamentalDomainStabInfty}.

            \begin{lemma}
                The map
                \begin{equation}
                    \begin{array}{ccccc}
                        \psi & : & \left( \mathbb{R}_+^{\ast} \right)^n & \longrightarrow & \left( \mathbb{R}_+^{\ast} \right)^n \\[0.5em]

                        && \left( y_1, \ldots, y_n \right) & \longmapsto & \left( y_1, \ldots, y_{n-1}, y_1 \cdot \ldots \cdot y_n \right)
                    \end{array}
                \end{equation}
                is a diffeomorphism which satisfies $\psi \left( B \right) \, = \, \left( \mathbb{R}_+^{\ast} \right)^{n-1} \times \left\{ 1 \right \}$, and whose Jacobian at $\left( y_1, \ldots, y_n \right)$ equals $y_1 \cdot \ldots \cdot y_{n-1}$. Furthermore, we have
                \begin{equation}
                    \begin{array}{lll}
                        \psi \left( \varepsilon \cdot \left( y_1, \ldots, y_n \right) \right) & = & \left( \varepsilon \cdot \left( y_1, \ldots, y_{n-1} \right), y_1 \cdot \ldots \cdot y_n \right).
                    \end{array}
                \end{equation}
            \end{lemma}

            \begin{proof}
                This is a direct computation.
            \end{proof}

            \begin{proposition}
                Denote by $G$ a fundamental domain for the lattice $\mathcal{O}_K^{\times,+}$ embedded in $\mathbb{R}^{n-1}$ as in proposition \ref{prop:OKTimesPlusLatticeRNMinus1}. The set
                \begin{equation}
                \label{eq:fundamentalDomainF}
                    \begin{array}{lll}
                        F & = & \psi^{-1} \left( \lambda^{-1} \left( G \right) \times \mathbb{R}_+^{\ast} \right)
                    \end{array}
                \end{equation}
                is a fundamental domain for the action of $\mathcal{O}_K^{\times,+}$ on $\left( \mathbb{R}_+^{\ast} \right)^n$. For any $r>0$, we further have
                \begin{equation}
                \label{eq:fundamentalDomainFr}
                    \begin{array}{lll}
                        F_r & = & \psi^{-1} \left( \lambda^{-1} \left( G \right) \times \left] 1/r^2, + \infty \right[ \right) \\[1em]

                        & = & \left \{ \left( y_1, \ldots, y_n \right) \in \left( \mathbb{R}_+^{\ast} \right)^n \; \middle \vert \; \left( y_1 \cdot \ldots \cdot y_n \right)^{-1/2} \, < \, r \right \}.
                    \end{array}
                \end{equation}
            \end{proposition}

            \begin{proof}
                This is a direct computation.
            \end{proof}

        \subsubsection{Distance to the cusps}

            The cusps can be thought of as ``points at infinity'' when considering the action of Hilbert modular groups, but there is a way to assign a meaningful distance between points in $\mathbb{H}^n$ and cusps.

            \begin{definition}
            \label{def:functionMu}
                Consider a cusp $c = \left[ \alpha : \beta \right] \in \mathbb{P}^1 \left( K \right)$, with $\alpha, \beta \in \mathcal{O}_K$. The function $\mu \left( \cdot, c \right)$ is defined as
                \begin{equation}
                    \begin{array}{ccccc}
                        \mu \left( \cdot, c \right) & : & \mathbb{H}^n & \longrightarrow & \mathbb{R}_+^{\ast} \\[1em]
                        && \tau & \longmapsto & \displaystyle \frac{N \left( \alpha, \beta \right)^2 N \left( \Imag \tau \right)}{\left \vert N \left( - \beta \tau + \alpha \right) \right \vert^2}
                    \end{array},
                \end{equation}
                where we have set $N \left( \alpha, \beta \right) \, = \, \# \left( \mathcal{O}_K / \left( \alpha \mathcal{O}_K + \beta \mathcal{O}_K \right) \right)$, and $N \left( w \right) \, = \, w_1 \ldots w_n$ for any $w \, = \, \left( w_1, \ldots, w_n \right) \in \mathbb{C}^n$. The function
                \begin{equation}
                    \begin{array}{ccccc}
                        \mu \left( \cdot, c \right)^{-1/2} & : & \mathbb{H}^n & \longrightarrow & \mathbb{R}_+^{\ast}
                    \end{array}
                \end{equation}
                then represents the \textit{distance to the cusp $c$}.
            \end{definition}

            \begin{remark}
            \label{rmk:functionMuInfty}
                In particular, we have
                \begin{equation}
                    \begin{array}{lllll}
                        \mu \left( \tau, \infty \right) & = & N \left( \Imag \tau \right) & = & y_1 \ldots y_n
                    \end{array}
                \end{equation}
                for any $\tau = \left( x_1 + iy_1, \ldots, x_n + i y_n \right) \in \mathbb{H}^n$.
            \end{remark}

            \begin{proposition}
            \label{prop:functionMuActionGamma}
                For any $\gamma \in \widehat{\Gamma}_K$, any $c \in \mathbb{P}^1 \left( K \right)$, and any $\tau \in \mathbb{H}^n$, we have
                \begin{equation}
                    \begin{array}{lll}
                        \mu \left( \gamma \cdot \tau, \gamma \cdot c \right) & = & \mu \left( \tau, c \right).
                    \end{array}
                \end{equation}
            \end{proposition}

            \begin{proof}
                This is a direct computation.
            \end{proof}

            \begin{definition}
                For any $c \in \mathbb{P}^1 \left( K \right)$ and any $r > 0$, we set
                \begin{equation}
                    \begin{array}{lll}
                        B \left( c, r \right) & = & \left \{ \tau \in \mathbb{H}^n \;\middle \vert \; \mu \left( \tau, c \right)^{-1/2} \, < \, r \right \},
                    \end{array}
                \end{equation}
                which can be thought of as the ball of radius $r$ centred at the cusp $c$.
            \end{definition}

            \begin{remark}
                When compactifying Hilbert modular varieties, the balls $B \left( c,r \right)$ provide a system of neighbourhoods for the cusps.
            \end{remark}

            \begin{proposition}
                For any $c \in \mathbb{P}^1 \left( K \right)$ and any $r>0$, the stabilizer $\widehat{\Gamma}_{K,c}$ of the cusp $c$ acts on $B \left( c,r \right)$.
            \end{proposition}

            \begin{proof}
                This is a direct consequence of proposition \ref{prop:functionMuActionGamma}.
            \end{proof}

            \begin{proposition}
            \label{prop:computationVolumeBall}
                For any $r>0$, we have
                \begin{equation}
                \label{eq:volumeBallRadiusRInfty}
                    \begin{array}{lll}
                        \vol \left( \widehat{\Gamma}_{K, \infty} \backslash B \left( \infty, r \right) \right) & = & \displaystyle \sqrt{\Delta_K} \cdot \frac{2^{n-1}}{[ \mathcal{O}_K^{\times,+} : \mathcal{O}_K^{\times,2}]} \cdot R_K \cdot r^2,
                    \end{array}
                \end{equation}
                the volume being computed with respect to the Poincaré metric.
            \end{proposition}

            \begin{proof}
                Fix $r>0$. We set
                \begin{equation}
                    \begin{array}{lll}
                        E_r & = & \left \{ \tau \in E \; \middle \vert \; \mu \left( \tau, \infty \right)^{-1/2} \, < \, r \right \},
                    \end{array}
                \end{equation}
                which is a fundamental domain for the action of $\widehat{\Gamma}_{K, \infty}$ on $B \left( \infty, r \right)$, having introduced $E$ in proposition \ref{prop:fundamentalDomainStabInfty}. Using the decomposition
                \begin{equation}
                    \begin{array}{lll}
                        E_r & = & T \; \times \; F_r,
                    \end{array}
                \end{equation}
                where $F_r$ was defined in \eqref{eq:fundamentalDomainFr}, we have
                \begin{equation}
                    \begin{array}{llll}
                        \vol \left( E_r \right) & = & \displaystyle \int_{E_r} \mathrm{d}\mu \left( \tau \right) \\[2em]

                        & = & \displaystyle \vol T \; \int_{F_r} \frac{\mathrm{d}y_1 \ldots \mathrm{d}y_n}{y_1^2 \ldots y_n^2} \\[2em]

                        & = & \displaystyle \sqrt{\Delta_K} \int_{\psi^{-1} \left( \lambda^{-1} \left( G \right) \times \left] 1/r^2, \infty \right[ \right)} \frac{\mathrm{d}y_1 \ldots \mathrm{d}y_n}{y_1^2 \ldots y_n^2} \\[2em]

                        & = & \displaystyle \sqrt{\Delta_K} \int_{\lambda^{-1} \left( G  \right)} \frac{\mathrm{d}y_1 \ldots \mathrm{d}y_{n-1}}{y_1 \ldots y_{n-1}} \int_{1/r^2}^{\infty} \frac{\mathrm{d}y_n}{y_n^2} \\[2em]

                        & = & \displaystyle \sqrt{\Delta_K} \; \vol G \; r^2 \\[1em]

                        & = & \displaystyle \sqrt{\Delta_K} \cdot \frac{2^{n-1}}{[ \mathcal{O}_K^{\times,+} : \mathcal{O}_K^{\times,2}]} \cdot R_K \cdot r^2, \\[2em]
                    \end{array}
                \end{equation}
                the volume of $G$ being given by \eqref{eq:volumeActionOKTimesPlusRNMinus1}.
            \end{proof}

            \begin{remark}
            \label{rmk:quadraticHomogeneityVolume}
                From \eqref{eq:volumeBallRadiusRInfty}, one deduces the following equality
                \begin{equation}
                    \begin{array}{lll}
                        \vol \left( \widehat{\Gamma}_{K, \infty} \backslash B \left( \infty, r \right) \right) & = & r^2 \; \vol \left( \widehat{\Gamma}_{K, \infty} \backslash B \left( \infty, 1 \right) \right)
                    \end{array}
                \end{equation}
                for any $r > 0$.
            \end{remark}

            Let us end this section with a way to describe particular fundamental domains for the action of $\widehat{\Gamma}_K$ on $\mathbb{H}^n$ using the function $\mu$.

            \begin{definition}
            \label{def:sphereInfluence}
                Let $c \in \mathbb{P}^1 \left( K \right)$. The \textit{sphere of influence} of $c$ is defined as
                \begin{equation}
                    \begin{array}{lll}
                        S_c & = & \left \{ \tau \in \mathbb{H}^n \; \middle \vert \; \forall c' \in \mathbb{P}^1 \left( K \right) \setminus \left \{ c \right \}, \; \mu \left( \tau, c \right) \, \geqslant \, \mu \left( \tau, c' \right) \right \}.
                    \end{array}
                \end{equation}
            \end{definition}

            \begin{proposition}
                Let $\gamma \in \widehat{\Gamma}_K$. We have
                \begin{equation}
                    \begin{array}{lll}
                        \gamma \cdot S_c^{\circ} \, \cap \, S_c^{\circ} \; \neq \; \emptyset & \iff & \gamma \in \widehat{\Gamma}_{K, c},
                    \end{array}
                \end{equation}
                where $S_c^{\circ}$ denotes the interior of $S_c$.
            \end{proposition}

            \begin{proof}
                This is explained in \cite[Proposition I.2.3]{vanDerGeer:hilbert-modular-surfaces}.
            \end{proof}

            \begin{proposition}
            \label{prop:sphereInfluenceFundamentalDomain}
                Consider a set $c_1, \ldots, c_{h_K} \in \mathbb{P}^1 \left( K \right)$ of representatives of the cusps modulo the action of $\widehat{\Gamma}_K$. The set
                \begin{equation}
                    \begin{array}{lll}
                        S & = & \bigsqcup\limits_{j=1}^{h_K} \; \widehat{\Gamma}_{K,c_j} \backslash S_{c_j}
                    \end{array}
                \end{equation}
                is a fundamental domain for the action of $\widehat{\Gamma}_K$ on $\mathbb{H}^n$.
            \end{proposition}

            \begin{proof}
                This is explained in \cite[Section I.2]{vanDerGeer:hilbert-modular-surfaces}.
            \end{proof}

\section{A Minkowski-type theorem}

    This final section will be devoted to associating a rigid adelic space to any point in $\mathbb{H}^n$, and to seeing how some of the natural invariants studied by Gaudron in~\cite{gaudron:rigid-adelic-spaces} appear in this picture. The notations from previous sections will be adopted here.

    \subsection{\texorpdfstring{From points in $\mathbb{H}$ to positive-definite matrices}{From points in H to positive-definite matrices}}
    \label{subsec:posdeftoH}

        The correspondence between points in the Poincaré upper-half plane $\mathbb{H}$ and $2 \times 2$ symmetric positive-definite matrices is well-established, and for instance used in \cite{siegel:volume-fundamental-domain}. Being crucial in this section, let us recall the construction.

        \begin{proposition}
            The map
            \begin{equation}
                \begin{array}{ccccc}
                    \varphi & : & \mathbb{H} \times \mathbb{R}_+^{\ast} & \longrightarrow & S_2^{++} \left( \mathbb{R} \right) \\[1em]
                    && \left( \tau, \, \lambda \right) & \longmapsto & \displaystyle \frac{\sqrt{\lambda}}{y} \left( \begin{array}{cc} x^2+y^2 & x \\[0.4em] x & 1 \end{array} \right)
                \end{array},
            \end{equation}
            where we write $\tau = x+iy$, is invertible, and its inverse is given by
            \begin{equation}
                \begin{array}{ccccccc}
                    \psi & : & \multicolumn{3}{c}{S_2^{++} \left( \mathbb{R} \right)} & \longrightarrow & \mathbb{H} \times \mathbb{R}_+^{\ast} \\[1em]

                    && S & = & \left( \begin{array}{cc} u & v \\[0.4em] v & w \end{array} \right) & \longmapsto & \displaystyle \left( \frac{v + i \sqrt{\det S}}{w}, \, \det S \right)
                \end{array}.
            \end{equation}
        \end{proposition}

        \begin{proof}
            The equalities $\varphi \circ \psi = \id_{S_2^{++} \left( \mathbb{R} \right)}$ and $\psi \circ \varphi = \id_{\mathbb{H} \times \mathbb{R}_+^{\ast}}$ can be checked by direct computations.
        \end{proof}

        \begin{definition}
            For any $\tau \in \mathbb{H}$, we set $S \left( \tau \right) = \varphi \left( \tau, 1 \right)$, and
            \begin{equation}
                \begin{array}{ccc}
                    T \left( \tau \right) & = & \displaystyle \frac{1}{\sqrt{y}} \left( \begin{array}{cc} y & 0 \\[0.4em] x & 1 \end{array} \right),
                \end{array}
            \end{equation}
            having once again written $\tau = x + iy$.
        \end{definition}

        \begin{proposition}
            For any $\tau \in \mathbb{H}$, we have
            \begin{equation}
                \begin{array}{lll}
                    S \left( \tau \right) & = & \transpose{T \left( \tau \right)} T \left( \tau \right).
                \end{array}
            \end{equation}
        \end{proposition}

        \begin{proof}
            This is a direct computation.
        \end{proof}

        \begin{proposition}
        \label{prop:phiGroupActions}
            For any matrix $P \in PSL_2 \left( \mathbb{R} \right)$, and any $\left( \tau, \lambda \right) \in \mathbb{H} \times \mathbb{R}_+^{\ast}$, we have
            \begin{equation}
                \begin{array}{ccc}
                    \varphi \left( P \cdot \tau, \, \lambda \right) & = & P \varphi \left( \tau, \, \lambda \right) \transpose{P}.
                \end{array}
            \end{equation}
        \end{proposition}

        \begin{proof}
            Consider a matrix
            \begin{equation}
                \begin{array}{lllll}
                    P & = & \left( \begin{array}{cc} a & b \\ c & d \end{array} \right) & \in & PSL_2 \left( \mathbb{R} \right).
                \end{array}
            \end{equation}
            Writing $\tau \, = \, x+iy$, we have
            \begin{equation}
                \begin{array}{lllll}
                    \displaystyle P \cdot \tau & = & \displaystyle \frac{a \tau + b}{c \tau + d} & = & \displaystyle \frac{\left( a \tau + b \right) \left( c \overline{\tau} + d \right)}{\left \vert c \tau + d \right \vert^2} \\[2em]
                    
                    &&& = & \displaystyle \frac{ac \left \vert \tau \right \vert^2 + ad \tau + bc \overline{\tau} + bd}{\left \vert c \tau + d \right \vert^2} \\[2em]

                    &&& = & \displaystyle \frac{ac \left( x^2 + y^2 \right) + \left( ad + bc \right) x + bd + iy}{\left( cx+d \right)^2 + c^2 y^2} \\[2em]

                    &&& = & \displaystyle \frac{ac \left( x^2 + y^2 \right) + \left( ad + bc \right) x + bd + iy}{c^2 \left( x^2+y^2 \right) + 2 dc x + d^2},
                \end{array}
            \end{equation}
            and we also have
            \begin{equation}
                \begin{array}{llll}
                    \multicolumn{2}{l}{P \varphi \left( \tau, \, \lambda \right) \transpose{P}} & = & \displaystyle \frac{\sqrt{\lambda}}{y} \left( \begin{array}{cc} a & b \\[0.4em] c & d \end{array} \right) \left( \begin{array}{cc} x^2+y^2 & x \\[0.4em] x & 1 \end{array} \right) \left( \begin{array}{cc} a & c \\[0.4em] b & d \end{array} \right) \\[2em]

                    = & \multicolumn{3}{l}{\displaystyle \frac{\sqrt{\lambda}}{y} \left( \begin{array}{cc} a \left( x^2+y^2 \right) +bx & ax+b \\[0.4em] c \left( x^2+y^2 \right)+dx & cx+d \end{array} \right) \left( \begin{array}{cc} a & c \\[0.4em] b & d \end{array} \right)} \\[2em]

                    = & \multicolumn{3}{l}{\displaystyle \frac{\sqrt{\lambda}}{y} \left( \begin{array}{cc} a^2 \left( x^2+y^2 \right) + 2abx + b^2 & \substack{\displaystyle ac \left( x^2+y^2 \right) + \left( ad+bc \right)x \\ \displaystyle + bd} \\[0.8em] \substack{\displaystyle ac \left( x^2+y^2 \right) + \left( ad+bc \right)x \\ \displaystyle + bd} & c^2 \left( x^2+y^2 \right) + 2dcx + d^2 \end{array} \right).}
                \end{array}
            \end{equation}
            We thus get
            \begin{equation}
                \begin{array}{lll}
                    \psi \left( P \varphi \left( \tau, \, \lambda \right) \transpose{P} \right) & = & \displaystyle \left( \frac{\frac{\sqrt{\lambda}}{y} \left( ac \left( x^2+y^2 \right) + \left( ad+bc \right) x + bd \right) + i \sqrt{\lambda}}{\frac{\sqrt{\lambda}}{y} \left( c^2 \left( x^2+y^2 \right) + 2dcx + d^2 \right)}, \; \lambda \right) \\[2em]

                    & = & \displaystyle \left( \frac{ac \left( x^2+y^2 \right) + \left( ad+bc \right) x + bd + iy}{c^2 \left( x^2+y^2 \right) + 2dcx + d^2}, \; \lambda \right) \\[2em]

                    & = & \left( P \cdot \tau, \; \lambda \right).
                \end{array}
            \end{equation}
        \end{proof}

        \begin{remark}
            As a direct consequence of proposition \ref{prop:phiGroupActions}, we have
            \begin{equation}
                \begin{array}{lll}
                    S \left( P \cdot \tau \right) & = & P S \left( \tau \right) \transpose{P}
                \end{array}
            \end{equation}
            for any $P \in PSL_2 \left( \mathbb{R} \right)$ and $\tau \in \mathbb{H}$.
        \end{remark}

    \subsection{\texorpdfstring{From points in $\mathbb{H}^n$ to rigid adelic spaces}{From points in Hn to rigid adelic spaces}}

        In the previous section, we saw how to obtain a symmetric, positive-definite matrix $S \left( \tau \right)$ from any point $\tau \in \mathbb{H}$. Applying this construction $n$ times thus naturally leads to an injective map
        \begin{equation}
            \begin{array}{ccc}
                \mathbb{H}^n & \longrightarrow & S_2^{++} \left( \mathbb{R} \right)^n.
            \end{array}
        \end{equation}
        Let us now make use of this map in order to construct rigid adelic spaces, as introduced in definition \ref{def:rigidAdelicSpaceDim2}. In the following, we consider a totally real number field $K$, and adopt the notations from previous sections.
        
        \begin{definition}
            For any point $\tau = \left( \tau_1, \ldots, \tau_n \right) \in \mathbb{H}^n$, we define the rigid adelic space $E_{\tau}$ as given by the matrices $\left( A_{\tau, v} \right)_{v \in V \left( K \right)}$, with
            \begin{equation}
                \begin{array}{lll}
                    A_{\tau, v} & = & \left \{ \begin{array}{ll} I_2 & \text{if } v \text{ is finite} \\[0.3em] T \left( \tau_j \right) & \text{if } v \text{ is infinite, associated with } \sigma_j \end{array} \right. .
                \end{array}
            \end{equation}
        \end{definition}

        \begin{proposition}
        \label{prop:HETauEquals1}
            For any $\tau \in \mathbb{H}^n$, we have
            \begin{equation}
                \begin{array}{lll}
                    H \left( E_{\tau} \right) & = & 1.
                \end{array}
            \end{equation}
        \end{proposition}

        \begin{proof}
            This is a consequence of the fact that we have
            \begin{equation}
                \begin{array}{lll}
                    \det A_{\tau, v} & = & 1
                \end{array}
            \end{equation}
            for any $\tau \in \mathbb{H}^n$ and any $v \in V \left( K \right)$.
        \end{proof}

    \subsection{The relation between heights and distances to cusps}

        Recall that, in definition \ref{def:HeightPoint}, we attached to any point $z$ in a rigid adelic space $E$ a number~$H_E \left( z \right)$ called the height of $z$. In this paragraph, we will see that the heights of points of~$E_{\tau}$ are closely related with the distances between~$\tau$ and the cusps of $\mathbb{H}^n$ with respect to $\widehat{\Gamma}_K$, whose set is naturally identified with~$\mathbb{P}^1 \left( K \right)$. From now on, we will make the following assumption.

        \begin{assumption}
        \label{assumption:class-number-1}
            The class number $h_K$ of $K$ equals $1$.
        \end{assumption}

        \begin{definition}
            We define the involution $\iota$ of $\mathbb{P}^1 \left( K \right)$ by
            \begin{equation}
                \begin{array}{ccccc}
                    \iota & : & \mathbb{P}^1 \left( K \right) & \longrightarrow & \mathbb{P}^1 \left( K \right) \\[0.5em]
                    && \left[ x,y \right] & \longmapsto & \left[ y : -x \right]
                \end{array}.
            \end{equation}
        \end{definition}

        \begin{proposition}
        \label{prop:heightAndMuFunction}
            For any $\tau \in \mathbb{H}^n$, we have
            \begin{equation}
            \label{eq:heightAndMuFunction}
                \begin{array}{ccc}
                    H_{E_{\tau}} \circ i & = & \mu \left( \tau, \, \cdot \right)^{-1/2n},
                \end{array}
            \end{equation}
            where $\mu$ is the function introduced in definition \ref{def:functionMu}.
        \end{proposition}

        \begin{proof}
            Consider $\tau \in \mathbb{H}^n$, as well as $\left[ \alpha: \beta \right] \in \mathbb{P}^1 \left( K \right)$. Without loss of generality, we may assume that we have $\alpha, \beta \in \mathcal{O}_K$. The class number $h_K$ being equal to $1$, by assumption \ref{assumption:class-number-1}, we can find $d \in \mathcal{O}_K$ such that we have
            \begin{equation}
            \label{eq:gcdOfxAndy}
                \begin{array}{ccc}
                    \alpha \, \mathcal{O}_K + \beta \, \mathcal{O}_K & = & d \, \mathcal{O}_K,
                \end{array}
            \end{equation}
            \textit{i.e.} $d$ is the gcd of $\alpha$ and $\beta$. For any finite place $v$ of $K$, we then have
            \begin{equation}
                \begin{array}{lll}
                    \left \vert d \right \vert_v & = & \max \left( \left \vert \alpha \right \vert_v, \, \left \vert \beta \right \vert_v \right).
                \end{array}
            \end{equation}
            The product formula (see \eqref{eq:productFormula}) applied to $d$ then yields
            \begin{equation}
                \begin{array}{lll}
                    \displaystyle \prod\limits_{v \in V_f \left( K \right)} \left \vert d \right \vert_v^{n_v / n} & = & \displaystyle \prod\limits_{j=1}^n \left \vert \sigma_j \left( d \right) \right \vert^{-1/n},
                \end{array}
            \end{equation}
            since we have $n_v \, = \, 1$ if $v$ is infinite, and we get
            \begin{equation}
                \begin{array}{lll}
                    H_{E_{\tau}} \left( \left[ \alpha : \beta \right] \right) & = & \displaystyle \prod\limits_{\substack{v \in V \left( K \right) \\ v \text{ finite}}} \left( \max \left( \left \vert \alpha \right \vert_v, \, \left \vert \beta \right \vert_v \right) \right)^{n_v/n} \cdot \prod\limits_{j=1}^n \left \Vert \left( \alpha, \beta \right) \right \Vert_{\tau,j}^{1/n} \\[2em]

                    & = & \displaystyle \prod\limits_{\substack{v \in V \left( K \right) \\ v \text{ finite}}} \left \vert d \right \vert_v^{n_v/n} \cdot \prod\limits_{j=1}^n \left \Vert \left( \alpha, \beta \right) \right \Vert_{\tau,j}^{1/n} \\[2em]

                    & = & \displaystyle \prod\limits_{j=1}^n \left \Vert \left( \frac{\alpha}{d},\frac{\beta}{d} \right) \right \Vert_{\tau,j}^{1/n}.
                \end{array}
            \end{equation}
            Using \eqref{eq:gcdOfxAndy}, we now let $u,w \in \mathcal{O}_K$ be such that we have $\alpha u + \beta w \, = \, d$, and set
            \begin{equation}
                \begin{array}{lllll}
                    \gamma & = & \left( \begin{array}{cc} w & -u \\[0.4em] \alpha/d & \beta/d \end{array} \right) & \in & PSL_2 \left( \mathcal{O}_K \right)
                \end{array}
            \end{equation}
            For any $j \in \llbracket 1, n \rrbracket$, we have
            \begin{equation}
                \begin{array}{lllll}
                    \displaystyle \left \Vert \left( \frac{\alpha}{d},\frac{\beta}{d} \right) \right \Vert_{\tau,j} & = & \displaystyle \frac{1}{\sqrt{y_j}} \left \vert \sigma_j \left( \frac{\alpha}{d} \right)  \tau_j + \sigma_j \left( \frac{\beta}{d} \right) \right \vert & = & \Imag \left( \sigma_j \left( \gamma \right) \tau_j \right)^{-1/2},
                \end{array}
            \end{equation}
            thereby giving
            \begin{equation}
                \begin{array}{lllll}
                    H_{E_{\tau}} \left( \left[ \alpha : \beta \right] \right) & = & \displaystyle \prod\limits_{j=1}^n \Imag \left( \sigma_j \left( \gamma \right) \tau_j \right)^{-1/2n} & = & \mu \left( \gamma \tau, \, \infty \right)^{-1/2n} \\[1em]

                    &&& = & \mu \left( \tau, \gamma^{-1} \infty \right)^{-1/2n} \\[1em]

                    &&& = & \mu \left( \tau, \iota \left( \left[ \alpha : \beta \right] \right) \right)^{-1/2n}.
                \end{array}
            \end{equation}
            We get formula \eqref{eq:heightAndMuFunction} by using the fact that $\iota$ is an involution of $\mathbb{P}^1 \left( K \right)$.
        \end{proof}

        Proposition \ref{prop:heightAndMuFunction} will now allow us to relate the Roy--Thunder minima of $E_{\tau}$ to the distances between $\tau$ and its two closest cusps.

        \begin{definition}
            For any $\tau \in \mathbb{H}^n$, we set
            \begin{equation}
            \label{eq:defMu1}
                \begin{array}[t]{lll}
                    \mu_1 \left( \tau \right) & = & \max\limits_{c \in \mathbb{P}^1 \left( K \right)} \mu \left( \tau, \, c \right),
                \end{array}
            \end{equation}
            and, denoting by $c_{\tau} \in \mathbb{P}^1 \left( K \right)$ a cusp realizing the maximum in \eqref{eq:defMu1}, we also set
            \begin{equation}
                \begin{array}[t]{lll}
                    \mu_2 \left( \tau \right) & = & \max\limits_{\substack{c \in \mathbb{P}^1 \left( K \right) \\ c \neq c_{\tau}}} \mu \left( \tau, \, c \right).
                \end{array}
            \end{equation}
        \end{definition}

        \begin{remark}
            Note that there are points in $\mathbb{H}^n$ such that we have $\mu_1 \left( \tau \right) \, = \, \mu_2 \left( \tau \right)$. These will be important later.
        \end{remark}

        \begin{proposition}
        \label{prop:Lambda1Lambda2Mu}
            For any $\tau \in \mathbb{H}^n$, we have
            \begin{equation}
            \label{eq:Lambda1Lambda2Mu}
                \begin{array}{lll}
                    \Lambda_1 \left( E_{\tau} \right) & = & \mu_1 \left( \tau \right)^{-1/2n}, \\[0.5em]

                    \Lambda_2 \left( E_{\tau} \right) & = & \mu_2 \left( \tau \right)^{-1/2n}.
                \end{array}
            \end{equation}
        \end{proposition}

        \begin{proof}
            Let $\tau \in \mathbb{H}^n$. For the first equality in \eqref{eq:Lambda1Lambda2Mu}, we have
            \begin{equation}
                \begin{array}{lllll}
                    \Lambda_1 \left( E_{\tau} \right) & = & \min\limits_{\left( \alpha, \beta \right) \in E_{\tau}} H_{E_{\tau}} \left( \left[ \alpha : \beta \right] \right) & = & \min\limits_{\left( \alpha, \beta \right) \in E_{\tau}} H_{E_{\tau}} \circ \iota \left( \left[ \alpha : \beta \right] \right) \\[1em]
                    
                    &&& = & \min\limits_{\left[ \alpha : \beta \right] \in \mathbb{P}^1 \left( K \right)} \left( \mu \left( \tau, \left[ \alpha : \beta \right] \right)^{-1/2n} \right) \\[1em]

                    &&& = & \mu_1 \left( \tau \right)^{-1/2n}.
                \end{array}
            \end{equation}
            Let us move on to the second equality in \eqref{eq:Lambda1Lambda2Mu}. To that effect, let us note that, for any $\left( \alpha_1, \beta_1 \right), \left( \alpha_2, \beta_2 \right) \in E_{\tau} \setminus \left \{ 0 \right \}$ we have
            \begin{equation}
                \begin{array}{lll}
                    \multicolumn{3}{l}{\dim_K \vect \left( \left( \alpha_1, \beta_1 \right), \left( \alpha_2, \beta_2 \right) \right) \; = \; 2} \\[0.5em]
                    
                    \qquad \qquad \qquad \qquad \qquad & \iff & \left[ \alpha_1 ; \beta_1 \right] \; \neq \; \left[ \alpha_2 : \beta_2 \right] \text{ in } \mathbb{P}^1 \left( K \right).
                \end{array}
            \end{equation}
            Thus we have
            \begin{equation}
                \begin{array}{lll}
                    \Lambda_2 \left( E_{\tau} \right) & = & \min \left\{ \max \left( H_{E_{\tau}} \left( \left[ \alpha_1 : \beta_1 \right] \right), H_{E_{\tau}} \left( \left[ \alpha_2 : \beta_2 \right] \right) \right) \; \middle \vert \right. \\
                    
                    && \qquad \qquad \qquad \qquad \qquad \quad \left. \left[ \alpha_1 : \beta_1 \right] \neq \left[ \alpha_2 : \beta_2 \right] \text{ in } \mathbb{P}^1 \left( K \right) \right \} \\[1em]

                    & = & \min \left\{ \max \left( H_{E_{\tau}} \circ \iota \left( \left[ \alpha_1 : \beta_1 \right] \right), H_{E_{\tau}} \circ \iota \left( \left[ \alpha_2 : \beta_2 \right] \right) \right) \; \middle \vert \right. \\
                    
                    && \qquad \qquad \qquad \qquad \qquad \quad \left. \left[ \alpha_1 : \beta_1 \right] \neq \left[ \alpha_2 : \beta_2 \right] \text{ in } \mathbb{P}^1 \left( K \right) \right \} \\[1em]

                    & = & \min \left\{ \max \left( \mu \left( \tau, \, \left[ \alpha_1 : \beta_1 \right] \right)^{-1/2n}, \, \mu \left( \tau, \, \left[ \alpha_2 : \beta_2 \right] \right)^{-1/2n} \right) \; \vert \right. \\

                    && \qquad \qquad \qquad \qquad \qquad \quad \left. \left[ \alpha_1 : \beta_1 \right] \neq \left[ \alpha_2 : \beta_2 \right] \text{ in } \mathbb{P}^1 \left( K \right) \right \} \\[1em]

                    & = & \mu_2 \left( \tau \right)^{-1/2n}.
                \end{array}
            \end{equation}
        \end{proof}

    \subsection{A Minkowski-type theorem}

        The identification made in proposition \ref{prop:Lambda1Lambda2Mu} between the Roy--Thunder minima of $E_{\tau}$ and the distances between $\tau$ and its closest cusps allows us to state a version of Minkowski's second theorem on the function~$\mu$.

        \begin{theorem}
        \label{thm:minkowski}
            For any $\tau \in \mathbb{H}^n$, we have
            \begin{equation}
                \begin{array}{lllll}
                    \displaystyle \frac{1}{c_{II}^{\Lambda} \left( 2, K \right)^{4n}} & \leqslant & \mu_1 \left( \tau \right) \mu_2 \left( \tau \right) & \leqslant & 1,
                \end{array}
            \end{equation}
            where $c_{II}^{\Lambda} \left( 2, K \right)$ is the constant introduced in definition \ref{def:HermiteConstant}.
        \end{theorem}

        \begin{proof}
            This is a direct consequence of proposition \ref{prop:MinkowskiRoyThunder}, together with proposition~\ref{prop:HETauEquals1} and \eqref{eq:Lambda1Lambda2Mu}.
        \end{proof}

        Let us now see some corollaries of this result.

        \begin{corollary}
        \label{cor:separationCusps}
            Let $c, c' \in \mathbb{P}^1 \left( K \right)$ be two cusps. For any $\tau \in \mathbb{H}^n$, we have
            \begin{equation}
                \begin{array}{lll}
                    \mu \left( \tau, c \right) \; \geqslant \; 1 \text{ and } \mu \left( \tau, c' \right) \; \geqslant \; 1 & \implies & c \; = \; c'.
                \end{array}
            \end{equation}
            In other words, a point in $\mathbb{H}^n$ cannot be at a distance less than $1$ from two different cusps. As a consequence, we have
            \begin{equation}
            \label{eq:inclusionBInfinity1IntoSInfinity}
                \begin{array}{lll}
                    B \left( \infty, 1 \right) & \subseteq & S_{\infty},
                \end{array}
            \end{equation}
            where we recall that $S_{\infty}$ is the sphere of influence of $\infty$, and $B \left( \infty, 1 \right)$ is the ball centred at $\infty$ of radius $1$.
        \end{corollary}

        \begin{proof}
            By contraposition, assuming we have $c \neq c'$, the inequality
            \begin{equation}
                \begin{array}{lll}
                    \mu_1 \left( \tau \right) \mu_2 \left( \tau \right) & \leqslant & 1
                \end{array}
            \end{equation}
            gives $\mu_2 \left( \tau \right) \, \leqslant \, 1$, which implies that we have $\mu \left( \tau, c \right) \, \leqslant \, 1$ or $\mu \left( \tau, c' \right) \, \leqslant \, 1$.
        \end{proof}

        \begin{remark}
            This is an effective version of \cite[Lemma I.2.1]{vanDerGeer:hilbert-modular-surfaces}, at least under assumption \ref{assumption:class-number-1}.
        \end{remark}

        \begin{corollary}
        \label{cor:lowerBoundMu1}
            Let $c \in \mathbb{P}^1 \left( K \right)$ be a cusp. For any $\tau \in \mathbb{H}^n$, we have
            \begin{equation}
            \label{eq:lowerBoundMu1}
                \begin{array}{lll}
                    \displaystyle \mu_1 \left( \tau \right) & \geqslant & \displaystyle \frac{1}{c_{II}^{\Lambda} \left( 2, K \right)^{2n}}.
                \end{array}
            \end{equation}
            In other words, a point in $\mathbb{H}^n$ is always at distance at most $c_{II}^{\Lambda} \left( 2, K \right)^n$ from a cusp. As a consequence, we have
            \begin{equation}
            \label{eq:inclusionSInfinityIntoBInfinityC2K}
                \begin{array}{lll}
                    S_{\infty} & \subseteq & B \left( \infty, c_{II}^{\Lambda} \left( 2, K \right)^n \right).
                \end{array}
            \end{equation}
        \end{corollary}

        \begin{proof}
            The inequality
            \begin{equation}
                \begin{array}{lllll}
                    \displaystyle \frac{1}{c_{II}^{\Lambda} \left( 2, K \right)^{4n}} & \leqslant & \mu_1 \left( \tau \right) \mu_2 \left( \tau \right) & \leqslant & \mu_1 \left( \tau \right)^2
                \end{array}
            \end{equation}
            directly gives the result.
        \end{proof}

        \begin{remark}
            This is an effective and uniform version of \cite[Lemma I.2.2]{vanDerGeer:hilbert-modular-surfaces}, at least under assumption \ref{assumption:class-number-1}. By definition of the constant $c_{II}^{\Lambda} \left( 2, K \right)$, the lower bound in~\eqref{eq:lowerBoundMu1} is optimal. Thus, using the proof of \cite[Lemma I.2.2]{vanDerGeer:hilbert-modular-surfaces}, we have
            \begin{equation}
                \begin{array}{lll}
                    \displaystyle \frac{1}{2^n \Delta_K} & \leqslant & \displaystyle \frac{1}{c_{II}^{\Lambda} \left( 2, K \right)^{2n}},
                \end{array}
            \end{equation}
            which gives
            \begin{equation}
                \begin{array}{lll}
                    c_{II}^{\Lambda} \left( 2, K \right) & \leqslant & \sqrt{2} \Delta_K^{1/2n}.
                \end{array}
            \end{equation}
            This is precisely the estimate on the constant $c_{II}^{\Lambda} \left( 2, K \right)$ given by Gaudron in \cite{gaudron:rigid-adelic-spaces}, and proved by Gaudron and Rémond in \cite[Proposition 5.1]{gaudron-remond:corps-siegel}, using methods which are not (directly) related to the methods used by van der Geer in \cite{vanDerGeer:hilbert-modular-surfaces}.
        \end{remark}

        \begin{proposition}
            The boundary $\partial S_{\infty}$ of the sphere of influence $S_{\infty}$ satisfies
            \begin{equation}
                \begin{array}{lll}
                    \partial S_{\infty} & \subset & B \left( \infty, c_{II}^{\Lambda} \left( 2, K \right)^n \right) \setminus B \left( \infty, 1 \right).
                \end{array}
            \end{equation}
        \end{proposition}

        \begin{proof}
            Recall that, in definition \ref{def:sphereInfluence}, we defined the sphere of influence $S_\infty$ as
            \begin{equation}
                \begin{array}{lll}
                    S_{\infty} & = & \left \{ \tau \in \mathbb{H}^n \; \middle \vert \; \mu \left( \tau, \infty \right) \, = \, \mu_1 \left( \tau \right) \right \},
                \end{array}
            \end{equation}
            \textit{i.e.} as the set of points in $\mathbb{H}^n$ which are closer to $\infty$ than to any other cusps. The boundary of $S_{\infty}$ is then the set of points in $S_{\infty}$ which are equidistant from $\infty$ and another cusp, which means that we have
            \begin{equation}
                \begin{array}{lll}
                    \partial S_{\infty} & = & \left \{ \tau \in \mathbb{H}^n \; \middle \vert \; \mu \left( \tau, \infty \right) \, = \, \mu_1 \left( \tau \right) \, = \, \mu_2 \left( \tau \right) \right \}.
                \end{array}
            \end{equation}
            A point $\tau \in S_{\infty}$ must then satisfy, by theorem \ref{thm:minkowski},
            \begin{equation}
            \label{eq:minkowskiBoundarySInf}
                \begin{array}{lllll}
                    \displaystyle \frac{1}{c_{II}^{\Lambda} \left( 2, K \right)^{4n}} & \leqslant & \mu_1 \left( \tau \right)^2 & \leqslant & 1,
                \end{array}
            \end{equation}
            which yields the result, after applying an exponent $1/4$ to each part of \eqref{eq:minkowskiBoundarySInf}.
        \end{proof}

    \subsection{An interesting class of integrals}

        Using the Minkowski-type theorem presented in theorem \ref{thm:minkowski}, we will now derive interesting lower- and upper-bounds for normalized integrals of the type
        \begin{equation}
        \label{eq:targetIntegrals}
            \begin{array}{c}
                \displaystyle \frac{1}{\vol \left( \widehat{\Gamma}_K \backslash \mathbb{H}^n \right)} \; \int_{\widehat{\Gamma}_K \backslash \mathbb{H}^n} \; \mu_1 \left( \tau \right)^t \; \mathrm{d}\mu \left( \tau \right),
            \end{array}
        \end{equation}
        where $t$ is a real number to be adjusted later, and $\widehat{\Gamma}_K$ is the extended Hilbert modular group. These integrals appear in \cite[Section 4]{frey-lefourn-lorenzo:height-estimates}.

        \begin{definition}
            The partial volume function is defined as
            \begin{equation}
                \begin{array}{ccccc}
                    g & : & \mathbb{R}_+^{\ast} & \longrightarrow & \mathbb{R}_+^{\ast} \\[0.5em]
                    && x & \longmapsto & \vol \left( \widehat{\Gamma}_{K, \infty} \backslash \left(  S_{\infty} \cap B \left( \infty, x \right) \right) \right)
                \end{array}
            \end{equation}
        \end{definition}

        \begin{proposition}
            The partial volume function $g$ is increasing, and satisfies
            \begin{equation}
                \begin{array}{lll}
                    g \left( x \right) & = & \left \{ \begin{array}{ll} \displaystyle  \sqrt{\Delta_K} \cdot \frac{2^{n-1}}{[ \mathcal{O}_K^{\times, +} : \mathcal{O}_K^{\times, 2} ]} \cdot R_K \cdot x^2 & \text{if } x\leqslant 1 \\[2em] \displaystyle \vol \left( \widehat{\Gamma}_K \backslash \mathbb{H}^n \right) & \text{if } x \geqslant c_{II}^{\Lambda} \left( 2, K \right)^n \end{array} \right.
                \end{array}
            \end{equation}
        \end{proposition}

        \begin{proof}
            The formula giving $g \left( x \right)$ whenever we have $x \leqslant 1$ is a consequence of \eqref{eq:inclusionBInfinity1IntoSInfinity}, which gives
            \begin{equation}
                \begin{array}{lll}
                    S_{\infty} \cap B \left( \infty, x \right) & = & B \left( \infty, x \right)
                \end{array}
            \end{equation}
            for any $x \leqslant 1$, and thus
            \begin{equation}
                \begin{array}{lllll}
                    g \left( x \right) & = & \vol \left( \widehat{\Gamma}_{K, \infty} \backslash B \left( \infty, x \right) \right) & = & \sqrt{\Delta_K} \cdot \frac{2^{n-1}}{[ \mathcal{O}_K^{\times, +} : \mathcal{O}_K^{\times, 2} ]} \cdot R_K \cdot x^2
                \end{array}
            \end{equation}
            for these real numbers $x$, as a consequence of proposition \ref{prop:computationVolumeBall}. Similarly, the formula which holds when we have $x \geqslant c_{II}^{\Lambda} \left( 2, K \right)^n$ is a consequence of \eqref{eq:inclusionSInfinityIntoBInfinityC2K}, which gives
            \begin{equation}
                \begin{array}{lll}
                    S_{\infty} \cap B \left( \infty, x \right) & = & S_{\infty}
                \end{array}
            \end{equation}
            for these values of $x$, and of proposition \ref{prop:sphereInfluenceFundamentalDomain}.
        \end{proof}

        \begin{remark}
            As a consequence of remark \ref{rmk:quadraticHomogeneityVolume}, note that we have
            \begin{equation}
                \begin{array}{lll}
                    g \left( x \right) & = & x^2 \, g \left( 1 \right)
                \end{array}
            \end{equation}
            for any $x \leqslant 1$.
        \end{remark}

        \begin{remark}
            Computing $g$ between $1$ and $c_{II}^{\Lambda} \left( 2, K \right)^n$ is quite complicated, especially without a precise description of the sphere of influence $S_{\infty}$. For the purposes of this paper, using monotonicity to get lower- and upper-bounds on this region will be enough.
        \end{remark}

        \begin{proposition}
        \label{prop:integralMu1ExpT}
            For any $0 \leqslant t < 1$, we have
            \begin{equation}
            \label{eq:integralMu1ExpT}
                \begin{array}{lll}
                    \multicolumn{3}{l}{\displaystyle \int_{\widehat{\Gamma}_K \backslash \mathbb{H}^n} \; \mu_1 \left( \tau \right)^t \; \mathrm{d}\mu \left( \tau \right)} \\[0.5em]
                    
                    \qquad \qquad \qquad & = & \displaystyle \frac{\vol \left( \widehat{\Gamma}_K \backslash \mathbb{H}^n \right)}{c_{II}^{\Lambda} \left( 2, K \right)^{2nt}} + 2t \int_{0}^{c_{II}^{\Lambda} \left( 2, K \right)^n} g \left( x \right) x^{-2t-1} \mathrm{d}x.
                \end{array}
            \end{equation}
        \end{proposition}

        \begin{proof}
            The idea behind this proof is to express the left-hand side of~\eqref{eq:integralMu1ExpT} as a Stieltjes integral (see \cite[Chapter 6]{rudin:maths-analysis}) with $g$ as the integrator. We have
            \begin{equation}
                \begin{array}[t]{lllll}
                    \displaystyle \int_{\widehat{\Gamma}_K \backslash \mathbb{H}^n} \; \mu_1 \left( \tau \right)^t \; \mathrm{d}\mu \left( \tau \right) & = & \displaystyle \lim\limits_{\varepsilon \rightarrow 0^+} \int_{\substack{\left \{ \tau \in \widehat{\Gamma}_K \backslash \mathbb{H}^n \; \vert \right. \\ \left. \mu_1 \left( \tau \right)^{-1/2} > \varepsilon \right \}}} \; \mu_1 \left( \tau \right)^t \; \mathrm{d}\mu \left( \tau \right).
                \end{array}
            \end{equation}
            Now, consider a real number $\varepsilon > 0$, and a partition $\left( T_k \right)_{k=1,\ldots, m}$ of $\left[ \varepsilon, c_{II}^{\Lambda} \left( 2, K \right)^n \right]$. Using corollary \ref{cor:lowerBoundMu1}, we have
            \begin{equation}
            \label{eq:breakIntegral}
                \begin{array}{lll}
                    \displaystyle \int_{\substack{\left \{ \tau \in \widehat{\Gamma}_K \backslash \mathbb{H}^n \; \vert \right. \\ \left. \mu_1 \left( \tau \right)^{-1/2} > \varepsilon \right \}}} \; \mu_1 \left( \tau \right)^t \; \mathrm{d}\mu \left( \tau \right) & \hspace{-0.6em} = \hspace{-0.6em} & \displaystyle \sum\limits_{k=1}^m \int_{\substack{\left \{ \tau \in \widehat{\Gamma}_K \backslash \mathbb{H}^n \; \vert \right. \\ \left. T_k < \mu_1 \left( \tau \right)^{-1/2} < T_{k+1} \right \}}} \; \mu_1 \left( \tau \right)^t \; \mathrm{d}\mu \left( \tau \right),
                \end{array}
            \end{equation}
            and we can control each integral in the right-hand side of \eqref{eq:breakIntegral}, as we have
            \begin{equation}
                \begin{array}{lllll}
                    \multicolumn{5}{l}{\displaystyle \sum\limits_{k=1}^m \frac{1}{T_{k+1}^{2t}} \left( g \left( T_{k+1} \right) - g \left( T_k \right) \right)} \\[2em]
                    
                    \qquad \qquad & < & \multicolumn{3}{l}{\displaystyle \sum\limits_{k=1}^m \int_{\substack{\left \{ \tau \in \widehat{\Gamma}_K \backslash \mathbb{H}^n \; \vert \right. \\ \left. T_k < \mu_1 \left( \tau \right)^{-1/2} < T_{k+1} \right \}}} \; \mu_1 \left( \tau \right)^t \; \mathrm{d}\mu \left( \tau \right)} \\[2em]

                    && \qquad \qquad  & < & \displaystyle \sum\limits_{k=1}^m \frac{1}{T_{k}^{2t}} \left( g \left( T_{k+1} \right) - g \left( T_k \right) \right),
                \end{array}
            \end{equation}
            and taking a limit on the partitions yields
            \begin{equation}
                \begin{array}{lll}
                    \displaystyle \sum\limits_{k=1}^m \int_{\substack{\left \{ \tau \in \widehat{\Gamma}_K \backslash \mathbb{H}^n \; \vert \right. \\ \left. T_k < \mu_1 \left( \tau \right)^{-1/2} < T_{k+1} \right \}}} \; \mu_1 \left( \tau \right)^t \; \mathrm{d}\mu \left( \tau \right) & = & \displaystyle \int_{\varepsilon}^{c_{II}^{\Lambda} \left( 2, K \right)^n} x^{-2t} \; \mathrm{d} g\left( x \right).
                \end{array}
            \end{equation}
            Note at this stage that $g$ may not be differentiable, so the last integral cannot directly be written as a Riemann integral. Using the integration by parts formula for Stieltjes integrals, which can be seen at the partition level as Abel's summation formula, and assuming we have $\varepsilon < 1$, we get
            \begin{equation}
                \begin{array}{lll}
                    \multicolumn{3}{l}{\displaystyle \int_{\varepsilon}^{c_{II}^{\Lambda} \left( 2, K \right)^n} x^{-2t} \; \mathrm{d} g\left( x \right)} \\[2em]
                    
                    & = & \displaystyle \frac{g \left( c_{II}^{\Lambda} \left( 2, K \right)^n \right)}{c_{II}^{\Lambda} \left( 2, K \right)^{2nt}} - \frac{g \left( \varepsilon \right)}{\varepsilon^{2t}} + 2t \int_{\varepsilon}^{c_{II}^{\Lambda} \left( 2, K \right)^n} g \left( x \right) x^{-2t-1} \mathrm{d}x \\[2em]

                    & = & \displaystyle \frac{\vol \left( \widehat{\Gamma}_K \backslash \mathbb{H}^n \right)}{c_{II}^{\Lambda} \left( 2, K \right)^{2nt}} - \varepsilon^{2 \left( 1-t \right)} g \left( 1 \right) + 2t \int_{\varepsilon}^{c_{II}^{\Lambda} \left( 2, K \right)^n} g \left( x \right) x^{-2t-1} \mathrm{d}x. \\[1.5em]
                \end{array}
            \end{equation}
            Having assumed that we have $0 \leqslant t < 1$, we get
            \begin{equation}
                \begin{array}{lll}
                    \lim\limits_{\varepsilon \rightarrow 0^+} \varepsilon^{2 \left( 1-t \right)} g \left( 1 \right) & = & 0,
                \end{array}
            \end{equation}
            and the function
            \begin{equation}
                \begin{array}{lll}
                    x & \longmapsto & g \left( x \right) x^{-2t-1}
                \end{array}
            \end{equation}
            is integrable at $0^+$ by the computation of $g$ for $x \leqslant 1$. Taking the limit as $\varepsilon$ goes to $0^+$ thus completes the proof.
        \end{proof}

        The integral on the right-hand side of \eqref{eq:integralMu1ExpT} can now be dealt with.

        \begin{proposition}
        \label{prop:inequalitiesIntegralG}
            For any $0 \leqslant t < 1$, we have
            \begin{equation}
                \begin{array}{lllll}
                    \displaystyle \frac{t}{1-t} g \left( 1 \right) & \leqslant & \multicolumn{3}{l}{\displaystyle 2t \int_{0}^{c_{II}^{\Lambda} \left( 2, K \right)^n} g \left( x \right) x^{-2t-1} \mathrm{d}x} \\[1em]
                    
                    &&& \leqslant & \displaystyle \frac{t}{1-t} g \left( 1 \right) + \vol \left( \widehat{\Gamma}_K \backslash \mathbb{H}^n \right) \left( 1 - \frac{1}{c_{II}^{\Lambda} \left( 2, K \right)^{2nt}} \right).
                \end{array}
            \end{equation}
        \end{proposition}

        \begin{proof}
            The function $g$ being explicitly computed on $\left[ 0,1 \right]$, we need only control it on the interval $\left[ 1, c_{II}^{\Lambda} \left( 2, K \right)^n \right]$, where we use positivity and monotonicity to find
            \begin{equation}
                \begin{array}{lllll}
                    0 & \leqslant & g \left( x \right) & \leqslant & \vol \left( \widehat{\Gamma}_K \backslash \mathbb{H}^n \right).
                \end{array}
            \end{equation}
            Applying this to the integral we want to study gives
            \begin{equation}
                \begin{array}{llll}
                    0 & \leqslant & \multicolumn{2}{l}{\displaystyle 2t \int_{0}^{c_{II}^{\Lambda} \left( 2, K \right)^n} g \left( x \right) x^{-2t-1} \mathrm{d}x} \\[1em]

                    && = & \displaystyle 2t \int_{0}^{1} g \left( x \right) x^{-2t-1} \mathrm{d}x + 2t \int_{1}^{c_{II}^{\Lambda} \left( 2, K \right)^n} g \left( x \right) x^{-2t-1} \mathrm{d}x \\[1em]

                    && \leqslant & \displaystyle 2t g\left( 1 \right) \int_0^1 x^{-2t+1} \mathrm{d}x + 2t \vol \left( \widehat{\Gamma}_K \backslash \mathbb{H}^n \right) \int_{1}^{c_{II}^{\Lambda} \left( 2, K \right)^n} x^{-2t-1} \mathrm{d}x \\[1em]

                    && \leqslant & \displaystyle \frac{t}{1-t} g \left( 1 \right) + \vol \left( \widehat{\Gamma}_K \backslash \mathbb{H}^n \right) \left( 1 - \frac{1}{c_{II}^{\Lambda} \left( 2, K \right)^{2nt}} \right).
                \end{array}
            \end{equation}
        \end{proof}

        Using these partial results, we can now state and prove the desired lower- and upper-bounds on the integrals \eqref{eq:targetIntegrals}.

        \begin{theorem}
            For any $0 \leqslant t < 1$, we have
            \begin{equation}
                \begin{array}{lllll}
                    \multicolumn{5}{l}{\displaystyle \frac{1}{c_{II}^{\Lambda} \left( 2, K \right)^{2nt}} + \frac{t}{1-t} \cdot \frac{1}{c_{II}^{\Lambda} \left( 2, K \right)^{2n}}} \\[2em]
                    
                    \qquad \qquad & \leqslant & \displaystyle \frac{1}{\vol \left( \widehat{\Gamma}_K \backslash \mathbb{H}^n \right)} \; \int_{\widehat{\Gamma}_K \backslash \mathbb{H}^n} \; \mu_1 \left( \tau \right)^t \; \mathrm{d}\mu \left( \tau \right) & \leqslant & \displaystyle \frac{1}{1-t}.
                \end{array}
            \end{equation}
        \end{theorem}

        \begin{proof}
            Consider $0 \leqslant t < 1$. Using proposition \ref{prop:integralMu1ExpT}, we have
            \begin{equation}
                \begin{array}{lll}
                    \multicolumn{3}{l}{\displaystyle \frac{1}{\vol \left( \widehat{\Gamma}_K \backslash \mathbb{H}^n \right)} \; \int_{\widehat{\Gamma}_K \backslash \mathbb{H}^n} \; \mu_1 \left( \tau \right)^t \; \mathrm{d}\mu \left( \tau \right)} \\[2em]
                    
                    \qquad \qquad & = & \displaystyle \frac{1}{c_{II}^{\Lambda} \left( 2, K \right)^{2nt}} + \frac{2t}{\vol \left( \widehat{\Gamma}_K \backslash \mathbb{H}^n \right)} \int_{0}^{c_{II}^{\Lambda} \left( 2, K \right)^n} g \left( x \right) x^{-2t-1} \mathrm{d}x,
                \end{array}
            \end{equation}
            which, after applying proposition \ref{prop:inequalitiesIntegralG}, yields
            \begin{equation}
            \label{eq:inequalityIntegral}
                \begin{array}{lllll}
                    \multicolumn{5}{l}{\displaystyle \frac{1}{c_{II}^{\Lambda} \left( 2, K \right)^{2nt}} + \frac{t}{1-t} \cdot \frac{g \left( 1 \right)}{\vol \left( \widehat{\Gamma}_K \backslash \mathbb{H}^n \right)}} \\[2em]

                    \qquad \qquad \qquad & \leqslant & \multicolumn{3}{l}{\displaystyle \frac{1}{\vol \left( \widehat{\Gamma}_K \backslash \mathbb{H}^n \right)} \; \int_{\widehat{\Gamma}_K \backslash \mathbb{H}^n} \; \mu_1 \left( \tau \right)^t \; \mathrm{d}\mu \left( \tau \right)} \\[2em]

                    && \qquad \qquad \qquad & \leqslant & \displaystyle 1 + \frac{t}{1-t} \cdot \frac{g \left( 1 \right)}{\vol \left( \widehat{\Gamma}_K \backslash \mathbb{H}^n \right)}.
                \end{array}
            \end{equation}
            To conclude, we need the following observation on the value $g\left( 1 \right)$. We have
            \begin{equation}
                \begin{array}{llll}
                    g \left( 1 \right) & = & \vol \left( \widehat{\Gamma}_{K, \infty} \backslash B \left( \infty, 1 \right) \right) \\[1em]
                    
                    & \leqslant & \vol \left( \widehat{\Gamma}_K \backslash \mathbb{H}^n \right) & \text{by proposition \ref{eq:inclusionBInfinity1IntoSInfinity}} \\[1em]

                    & \leqslant & \vol \left( \widehat{\Gamma}_{K, \infty} \backslash B \left( \infty, c_{II}^{\Lambda} \left( 2, K \right)^n \right) \right) & \text{by proposition \ref{eq:inclusionSInfinityIntoBInfinityC2K}} \\[1em]

                    & \leqslant & c_{II}^{\Lambda} \left( 2, K \right)^{2n} g \left( 1 \right) & \text{by remark \ref{rmk:quadraticHomogeneityVolume}},
                \end{array}
            \end{equation}
            which readily gives
            \begin{equation}
            \label{eq:InequalityG1}
                \begin{array}{lllll}
                    \displaystyle \frac{1}{c_{II}^{\Lambda} \left( 2, K \right)^{2n}} & \leqslant & \displaystyle \frac{g \left( 1 \right)}{\vol \left( \widehat{\Gamma}_K \backslash \mathbb{H}^n \right)} & \leqslant & 1.
                \end{array}
            \end{equation}
            Putting \eqref{eq:InequalityG1} into \eqref{eq:inequalityIntegral} completes the proof.
        \end{proof}


\bibliographystyle{amsplain}

\bibliography{bibliography}

\end{document}